\documentclass{amsart}
\usepackage{amssymb}
\usepackage{hyperref}
\newcommand{\hnl}{\htmladdnormallink}
\renewcommand{\eqref}[1]{{\rm(\ref{#1})}}
\newcommand{\Ref}[1]{{\rm\ref{#1}}}
\newcommand{\cal}{\mathcal}
\newcommand{\bA}{{\mathbf A}}
\newcommand{\bB}{{\mathbf B}}
\newcommand{\bbC}{{\mathbb C}}
\newcommand{\bbR}{{\mathbb R}}
\newcommand{\bbZ}{{\mathbb Z}}

\newcommand{\bH}{{\mathbf H}}
\newcommand{\bI}{{\mathbf I}}

\newcommand{\bQ}{{\mathbf Q}}

\newcommand{\bU}{{\mathbf U}}
\newcommand{\bV}{{\mathbf V}}

\newcommand{\cA}{{\mathcal A}}
\newcommand{\cB}{{\mathcal B}}

\newcommand{\cG}{{\mathcal G}}

\newcommand{\cH}{{\mathcal H}}

\newcommand{\cK}{{\mathcal K}}

\newcommand{\cN}{{\mathcal N}}
\newcommand{\cO}{{\mathcal O}}

\newcommand{\cS}{{\mathcal S}}

\newcommand{\ri}{{\rm i}}

\newcommand{\sH}{{\sf H}}

\newcommand{\sQ}{{\sf Q}}
\newcommand{\no}{\nonumber}

\newcommand{\tA}{\widetilde{A}}

\newcommand{\tH}{{\widetilde{H}}}

\newcommand{\III}{{|}\mkern -1.6mu{|}\mkern -1.6mu{|}}

\newcommand{\Max}{\mathop{{\rm max}}\limits}
\newcommand{\modulo}{\mathop{\rm mod\,}}

\newcommand{\Int}{\displaystyle\int\limits}
\newcommand{\Sup}{\mathop{\rm sup}}

\newcommand{\reduction}[2]{#1 \biggr|_{#2}}

\newcommand{\Lim}{\mathop{{\rm lim}}\limits}
\newcommand{\diag}{\mathop{\rm diag}}
\newcommand{\dist}{\mathop{\rm dist}}
\newcommand{\Img}{\mathop{\rm Im}}

\newcommand{\lal}{{\langle}}
\newcommand{\ral}{{\rangle}}

\newcommand{\be}{\begin{equation}}
\newcommand{\ee}{\end{equation}}

\DeclareMathOperator{\tr}{tr}

\DeclareMathOperator{\spec}{spec}

\DeclareMathOperator{\ran}{ran}

\DeclareMathOperator{\dom}{dom}

\DeclareMathOperator*{\nlim}{{\mathit n}-lim}
\DeclareMathOperator*{\slim}{{\mathit s}-lim}
\DeclareMathOperator*{\wlim}{{\mathit w}-lim}
\allowdisplaybreaks
\numberwithin{equation}{section}

\newtheorem{theorem}{Theorem}[section]
\newtheorem{lemma}[theorem]{Lemma}
\newtheorem{corollary}[theorem]{Corollary}
\newtheorem{remark}[theorem]{Remark}
\newtheorem{hypothesis}[theorem]{Hypothesis}
\newtheorem{definition}[theorem]{Definition}
\newtheorem{example}[theorem]{Example}
\theoremstyle{remark}
\begin{document}
\title[Graph Subspaces and the Spectral Shift Function]
{Graph Subspaces and the Spectral Shift Function}

\author[Albeverio, Makarov, Motovilov]{Sergio Albeverio,
Konstantin A.~Makarov,  and  Alexander~K.~Motovilov  }

\address{\hnl{Sergio Albeverio}
{http://wiener.iam.uni-bonn.de/albeverio/albeverio.html},
Institut f\"ur Angewandte Mathematik,
Universit\"at Bonn, Wegelerstra{\ss}e 6, D-53115 Bonn, Germany}
\email{\hnl{albeverio@uni-bonn.de}{mailto:albeverio@uni-bonn.de} \newline
\indent{\it URL}
\hnl{\tt http://wiener.iam.uni-bonn.de/albeverio/albeverio.html}
{http://wiener.iam.uni-bonn.de/albeverio/albeverio.html}}

\address{\hnl{Konstantin A. Makarov}
{http://www.math.missouri.edu/people/kmakarov.html},
Department of Mathematics, University of
Missouri, Co\-lum\-bia, MO
65211, USA}
\email{\hnl{makarov@math.missouri.edu}
{mailto:makarov@math.missouri.edu} \newline
\indent{\it URL}
\hnl{\tt http://www.math.missouri.edu/people/kmakarov.html}
{http://www.math.missouri.edu/people/kmakarov.html}}

\address{\hnl{Alexander K. Motovilov}
{http://www.jinr.ru/~motovilv}, Bogoliubov Laboratory of
Theoretical Physics, JINR, Joliot-Curie str. 6, 141980 Dubna, Russia}
\email{\hnl{motovilv@thsun1.jinr.ru}
{mailto:motovilv@thsun1.jinr.ru} \newline
\indent{\it URL} \hnl{\tt http://www.jinr.ru/\~{}motovilv}
{http://www.jinr.ru/~motovilv}}


\begin{abstract}
We extend the concept of Lifshits--Krein spectral shift
function associated with a pair of self-adjoint operators to the
case of  pairs of (admissible) operators that are similar to
self-adjoint operators.  An operator $H$  is called admissible if:
(i)~there is a bounded operator $V$ with a bounded inverse such
that $ H=V^{-1}\widehat HV$ for some self-adjoint operator
$\widehat H$; (ii)~the operators $H$ and $\widehat H$ are
resolvent comparable, i.\,e., the difference of the resolvents of
$H$ and  $\widehat H$ is a trace class operator (for non-real
values of the spectral parameter); (iii)~$\tr(VR-RV)=0$ whenever
$R$ is bounded and the commutator $VR-RV$ is a trace class
operator.  The spectral shift function $\xi(\lambda,H,A)$
associated with the pair of resolvent comparable admissible
operators $(H,A)$ is introduced then by the equality
$\xi(\lambda,H,A)=\xi(\lambda,\widehat H,\widehat A) $ where
$\xi(\lambda,\widehat H,\widehat A)$ denotes the Lifshits--Krein
spectral shift function associated with the pair $(\widehat
H,\widehat A)$ of self-adjoint operators.  Our main result is
the following. Let ${\cH}_0$ and ${\cH}_1$ be separable Hilbert
spaces, $A_0$ a self-adjoint operator in ${\cal H}_0$, $A_1$ a
self-adjoint operator in ${\cal H}_1$, and $B_{ij}$ a bounded
operator from ${\cal H}_j$ to ${\cal H}_i$, $i=0,1$, $j=1-i$,
and $B_{10}=B_{01}^*$.  Assume that
$
\bQ=\left(\begin{array}{lr}
  0             &   Q_{01}  \\
  Q_{10}        &   0
\end{array}\right)
$
is a Hilbert-Schmidt operator in $\cH$ with $Q_{01}=-Q_{10}^*$, and
the block operator matrix
$$
\bH=\left(\begin{array}{lr}
  A_0             &   B_{01}  \\
  B_{10}        &   A_1
\end{array}\right)=
\left(\begin{array}{lr}
  A_0             &  0  \\
 0        &   A_1
\end{array}\right)+
\left(\begin{array}{lr}
  0            &  B_{01}  \\
 B_{10}       &   0
\end{array}\right)=\bA+\bB
$$
has reducing graph subspaces of the form
$
\{x_i\oplus Q_{ji}x_i:\,x_i\in {\cH}_i\}, \quad
i=0,1,\quad j=1-i.
$
If both $(\bH-z\bI)^{-1}-(\bA-z\bI)^{-1}$ and
$\bB\bQ(\bA-z\bI)^{-1}$  are trace class operators in $\cH$ for
some $z$,  $\Img(z)>0$, then  the operators $A_i+B_{ij}Q_{ji}$
 and $A_i$, $i=0,1,$ $j=1-i$, acting in the spaces $\cH_i$ are
resolvent comparable admissible operators. Moreover, the
spectral shift function associated with the pair $(\bH,\bA)$ of
self-adjoint block operator matrices $\bH$ and $\bA$ admits the
representation $\xi(\lambda,\bH,\bA)
=\xi(\lambda,A_0+B_{01}Q_{10},A_0)+
\xi(\lambda,A_1+B_{10}Q_{01},A_1).$  We point out the intervals
on which the spectral shift functions
$\xi(\lambda,A_i+B_{ij}Q_{ji},A_i), $ $i=0,1,$ $j=1-i$, vanish.
We also obtain new representations for the solution to the
operator Sylvester equation in the form of Stieltjes operator integrals
 and formulate sufficient criterion
 for the strong solvability of the operator Riccati equation that ensures
the existence of reducing graph subspaces for block operator matrices.
\bigskip

\noindent\normalsize \bf LANL e-print: \hnl{\bf math.SP/0105142}
{http://www.arXiv.org/abs/math.SP/0105142}


\end{abstract}

\maketitle

\newpage
\tableofcontents
\newpage
\section{Introduction}
\label{SIntro}
The spectral analysis of operator block matrices is an important
issue in operator theory and mathematical physics. The search
for invariant subspaces, the problem of block diagonalization,
the analytic continuation of the compressed resolvents into
unphysical sheets of the spectral parameter plane as well as the
study of trace formulas attracted considerable attention in the
past due to numerous applications to various problems of quantum
mechanics, magnetohydrodynamics, and areas of mathematical
physics (see \cite{AdL}, \cite{AdLMSr}, \cite{Goedbloed},
\cite{Li89}, \cite{MM89}, \cite{MenShk}, \cite{MotRem} and
references cited therein).

In this work we restrict ourselves to the study of self-adjoint
operator block matrices of the form
\begin{equation}
\label{bH}
{\bH}=\left(\begin{array}{lr}
  A_0        &   B_{01}  \\
  B_{10}     &   A_{1}
\end{array}\right)
\end{equation}
acting in the orthogonal sum  ${\cH}={\cH}_0\oplus{\cH}_1$ of
separable Hilbert spaces ${\cH}_0$ and ${\cH}_1$.  The entries
$A_i$, $i=0,1,$ are assumed to be self-adjoint operators in
$\cH_i$ on domains $\dom(A_i)$. The off-diagonal elements
$B_{ij}:\,\cH_j\to\cH_i$,  $i=0,1$,
$j=1-i$, $B_{01}=B_{10}^*$, are assumed to be  bounded operators.

Under these assumptions the
matrix $\bH$ is a self-adjoint operator in $\cH$ on
$\dom(\bH)=\dom(A_0)\oplus\dom(A_1)=\dom(\bA)$ where
$\bA=\diag\{A_0,A_1\}$.
We also use the notation
$$
\bH=\bA+\bB\quad \text{ where }\quad\bB=\left(\begin{array}{lr}
    0        &   B_{01}  \\
  B_{10}     &   0
\end{array}\right).
$$

In the circle of ideas concerning the block diagonalization
problem for block operator matrices \eqref{bH} the existence of
invariant graph subspaces plays a crucial role.
Recall that a subspace $\cG_i$, $i=0$ or $i=1$, is said to be a
graph subspace of $\cH$ associated with the decomposition
$\cH=\cH_0\oplus\cH_1$ if it is the graph of a (bounded)
operator $Q_{ji}$, $j=1-i$, mapping $\cH_i$ to $\cH_j$.

The existence of a reducing graph subspace for a block operator matrix
\eqref{bH} is equivalent to the existence of a bounded
off-diagonal strong solution $\bQ$
to the  operator Riccati equation
\begin{equation}
\label{riceq0}
\bQ\bA-\bA\bQ+\bQ\bB\bQ=\bB
\end{equation}
having the form
\begin{equation}
\label{QQ-Vved}
\bQ=\left(\begin{array}{lr}
0      &   Q_{01}\\
Q_{10} &   0
\end{array}\right),
\qquad Q_{10}=-Q_{01}^*.
\end{equation}

Given a strong solution \eqref{QQ-Vved} to the equation
\eqref{riceq0}, the operator matrix $\bH=\bA+\bB$  has invariant
graph subspaces
$\cG_0=\{x\in \cH\, :\,
P_{\cH_1}x=Q_{10}P_{\cH_0}x\}$ and $\cG_1=\{x\in \cH\, :\,
P_{\cH_0}x=Q_{01}P_{\cH_1}x\}$ where $P_{\cH_i}$ denote the
orthogonal projections in $\cH=\cH_0\oplus \cH_1$ onto the
channel subspaces $\cH_i$, $i=0,1$. As a consequence, ${\bf H}$
can be block diagonalized
$$
(\bI+\bQ)^{-1}\bH(\bI+\bQ)=\bA+\bB\bQ=
\left(\begin{array}{cc}
A_0+B_{01}Q_{10}  &   0\\
0 &   A_1+B_{10}Q_{01}
\end{array}\right)
$$
by the similarity transformation generated by the operator $\bI +\bQ$.
Under these circumstances the block-diagonalization problem for
$\bH$ by a unitary transformation  admits an ``explicit''
solution,
\begin{equation}
\label{block}
\bU^*\bH\bU=\left(\begin{array}{cc}
H_0      &   0\\
0 &   H_1
\end{array}\right),
\end{equation}
where $\bU$ is the unitary operator from the polar decomposition
$\bI+\bQ=\bU|\bI+\bQ|$, and the diagonal entries $H_i$, $i=0,1,$
are self-adjoint operators similar to $A_0+B_{01}Q_{10}$ and
$A_1+B_{10}Q_{01}$, respectively.

Therefore, typical problems of {\it qualitative}
perturbation theory such as the existence of the graph invariant
subspaces, as well as a possibility of the block diagonalization
can be reduced to  purely analytic questions  concerning the
solvability of operator Riccati equations. Extensive
bibliography is devoted to the subject.  Not pretending to be
complete we refer to \cite{AdL00}, \cite{AdL}, \cite{AdLMSr},
\cite{AdLT}, \cite{A83}, \cite{AMSau}, \cite{ALMS},
\cite{BDMc83}, \cite{BR}, \cite{D53}, \cite{DK69}, \cite{DK70},
\cite{DR74}, \cite{LR59}, \cite{MM81}, \cite{MM99},
\cite{MennMotTMF}, \cite{MenShk}, \cite{MotSPbWorkshop},
\cite{MotRem}, \cite{P91}, \cite{R56}.

An intriguing problem of {\it quantitative} perturbation theory
is the study of the relationship between geometrical
characteristics of rotations of the invariant subspaces and the
accompanying shifts of the spectrum under a given perturbation.
It is the development of the quantitative perturbation theory
for self-adjoint block operator matrices that is the main goal
of the present paper.

In this context, the most important numerical quantitative
spectral characteristics  is the
Lifshits-Krein spectral shift function
\cite{Li52}, \cite{Li56}, \cite{Kr83}, \cite{Kr89}, \cite{Kr53}.
Detailed reviews of results on the spectral shift function and
its applications were published by Birman and Yafaev
\cite{BY93a}, \cite{BY93}, \cite{Ya92} and by Birman and
Pushnitskii \cite{BP98}. For many more references the interested
reader can consult \cite{GS96}, \cite{GM98}, \cite{GM00},
\cite{GMN}, \cite{Pu98}, \cite{Pu99a}.  For recent results we
refer to \cite{GMM}, \cite{K00}, \cite{Pu99}, \cite{Pu00},
and \cite{Si98}.

The spectral shift function $\xi(\lambda;\bH,\bA)$ associated
with the pair $(\bH,\bA)$ of self-adjoint operators is usually
introduced by the trace formula
\begin{equation}
\label{trfor}
\tr\big (\varphi(\bH)-\varphi(\bA)\big )=\int_\bbR
d\lambda\,
\varphi'(\lambda)\xi(\lambda; \bH, \bA).
\end{equation}
The trace formula \eqref{trfor} holds for a rather extensive
class of functions $\varphi:\,\bbR\to\bbC$, including the class
$C_0^\infty(\bbR)$ of infinitely differentiable functions with a
compact support, provided that the self-adjoint operators $\bH$
and $\bA$ are resolvent comparable, that is, the difference of
their resolvents is a trace class operator.

In case of the block operator matrices the quantitative spectral
analysis outlined above has a series of specific features. In
particular, if the matrix $\bH$ admits a block diagonalization
as in \eqref{block}, one might expect the validity of the
following splitting representation for the spectral shift
function
\begin{equation}
\label{split}
\xi(\lambda; \bH, \bA)=\xi(\lambda; H_0, A_0)+
\xi(\lambda; H_1, A_1).
\end{equation}
However, a certain difficulty in this way is that the spectral
shift function associated with a pair of self-adjoint operators
is not {\it stable} with respect to unitary transformations of
its operator arguments.  That is, if $\bU$ is a unitary operator,
the representation
\begin{equation}
\label{xiUn}
\xi(\lambda; \bU^*\bH\bU, \bA)=
\xi(\lambda; \bH, \bA)
\end{equation}
fails to hold in general, even if both terms in (\Ref{xiUn}) are
well-defined (see Example \Ref{example}).

One of the main goals of the present paper is to extend the
concept of the spectral shift function to pairs of {\it
admissible} (similar to self-adjoint) operators (see Definition
\Ref{defadmis}) followed by  the proof of the splitting formula
\eqref{split} as well as the proof of its ``non-self-adjoint''
version
\begin{equation}
\label{Int-split}
\xi(\lambda,\bH,\bA)
=\xi(\lambda,A_0+B_{01}Q_{10},A_0)+\xi(\lambda,A_1+B_{10}Q_{01},A_1)
\end{equation}
in the Hilbert-Schmidt class perturbation theory.

It is worth mentioning that the splitting formula
\eqref{Int-split} connects a purely spectral
characteristics of the perturbation, the spectral shift function
$\xi(\lambda,\bH,\bA)$, with the geometry of mutual disposition of
the invariant graph subspaces of the operator matrix $\bH$
determined by the angular operator $\bQ$ (provided that the
reducing graph subspaces for $\bH$ exist).

The plan of the paper is as follows.

In Section \Ref{SecSylv} we compare different representations
for the solutions of the operator Sylvester equation
\eqref{syl} and obtain new representations for its strong
solution based on the operator Stieltjes integrals approach.
These are the representations \eqref{EqX} and \eqref{EqXZ}.

In Section \Ref{SecRic} we extend our key result of Section
\Ref{SecSylv} (Theorem \Ref{SylUnb}) to the case  of the
operator Riccati equation
\begin{equation}
\label{riri}
QA-CQ+QBQ=D
\end{equation}
with self-ajoint (possibly unbounded) $A$ and $C$ and bounded
$B$ and $D$.  One of our main results (see Theorem \Ref{QsolvN})
provides a series of new sufficient conditions that imply the
weak or strong solvability of \eqref{riri}. We
prove, in particular,  that if the operators $A$ and $C$ are
bounded and
$$
\sqrt{\|B\|\|D\|}<\frac{1}{\pi}\,\dist\{\spec(A),\spec(C)\},
$$
then \eqref{riri} has even an operator solution. This result is
optimal in the following sense:  in case where $D=B^*$ the best
possible constant $c$ in the inequality
$$
\|B\|<c\,\dist\{\spec(A),\spec(C)\}
$$
that implies the solvability of \eqref{riri} lies within the
interval $\left[\frac{1}{\pi},\sqrt{2}\right]$ (see Remark
\Ref{bound}).

In Section \Ref{SecSSF} we introduce the concept of a spectral
shift function for the pairs of admissible operators which are
similar to self-adjoint (see Definitions \Ref{defadmis} and
\Ref{defSSF1}). We relate our general concept of the spectral
shift function associated with pairs of operators similar to
self-adjoint to the one based on  the perturbation determinant
approach originally suggested by Adamjan and Langer in  the case
of trace class perturbations \cite{AdL00}.

In Section \Ref{SecGPP} we discuss invariant graph subspaces for
block operator matrices and link their existence with the
existence of strong solutions to the corresponding Riccati
equations (Lemma \Ref{princip} and Theorem \Ref{th42}).

In Section \Ref{SecKey}, under rather general assumptions  we
prove the splitting formulas \eqref{split} and  \eqref{Int-split}
(Theorem \Ref{th42}).

Section \Ref{SecSuff} is devoted to a detailed study of the case
where the spectra of the diagonal entries  $A_0$ and $A_1$ of
the operator matrix $\bH$ are separated. Based on  the results
of Sec. \Ref{SecRic}  we prove one of the central results of the
present paper (Theorem \Ref{trsch} and Corollary \Ref{final})
concerning the validity of the splitting formulas \eqref{split},
\eqref{Int-split} in case of  Hilbert-Schmidt perturbations
$\bB$: if the perturbation $\bB$ is sufficiently small in a
certain sense (see Hypotheses \Ref{HEnorm} and \Ref{HBpi}) and
the operators $\bH=\bA+\bB$ and $\bA$ are resolvent comparable,
then

\begin{enumerate}

\item[(i)] the  splitting formulas \eqref{split} and
\eqref{Int-split} hold;

\item[(ii)]  the following equalities are valid
\begin{align*}
\xi(\lambda; H_0, A_0)&=\xi(\lambda; A_0+B_{01}Q_{10}, A_0)=0,\quad
\,\,\text{ for a.\,e. }
\lambda\in \spec(A_{0})
\\
\xi(\lambda; H_1, A_1)&=\xi(\lambda; A_1+B_{10}Q_{10}, A_1)=0,\quad
\text{ for a.\,e. }
\lambda\in \spec(A_{1}).
\end{align*}

\end{enumerate}

\section{Sylvester equation}
\label{SecSylv}

The principal purpose of this section is to introduce a new
representation for the solution $X$ of the operator Sylvester equation
$$
XA-CX=Y.
$$
We also discuss and compare the known representation theorems for
solution $X$. For a detail exposition and introduction to the
subject we refer to the papers \cite{BDMc83}, \cite{BR},
\cite{DR74}, \cite{LR59}, \cite{P91}, \cite{R56} and references
therein.

In the following $\cB(\cH,\cK)$  denotes the Banach
space of linear bounded operators between Hilbert spaces
$\cH$ and $\cK$. By $\cB_p(\cH,\cK)$, $p\geq1$, we understand
the standard Schatten\,--\,von Neumann ideals of $\cB(\cH,\cK)$.
For $\cB(\cH,\cH)$ and $\cB_p(\cH,\cH)$ we use the corresponding
shorten notation $\cB(\cH)$ and $\cB_p(\cH)$. The
$\cB_p(\cH,\cK)$--norm of a bounded operator $T$ acting
from  $\cH$ to $\cK$ is denoted by $\|T\|_p$.

Given two Hilbert spaces $\cH$ and $\cK$, recall the concept
of symmetric normed ideals of  $\cB(\cH, \cK)$, following
\cite{GK}.
\begin{definition}\label{enorm}
A two-sided ideal  $\cS\subset\cB(\cH, \cK)$ is called a
symmetric normed ideal of $\cB(\cH, \cK)$ if it is closed with
respect to a norm ${\III}\cdot {\III}$ on $\cS$ which has the
following properties:

\begin{enumerate}
\item[{\rm(i)}] if $T\in \cS$, $K\in \cB(K), $ $H\in \cB(\cH)$, then
$KTH \in\cB(\cH, \cK)$ and ${\III}KTH{\III}\le \|K\|\, {\III}T{\III} \, \|H\|$;

\item[{\rm(ii)}] if $T$ is rank one then ${\III}T{\III}=\|T\|$.
\end{enumerate}

\end{definition}

For technical reasons we also assume that
\begin{enumerate}
\item[{\rm(iii)}] {\it if $T_n\in \cS$ with $\sup_n {\III}T_n{\III} <
 \infty$, and if $T_n \to A$ in the weak operator topology, then
 $A \in \cS$ and ${\III}A{\III}\le \sup_n {\III}T_n{\III}$.}
\end{enumerate}

Recall that if $\cK=\cH$ then for any symmetric normed
ideal $\cS$ possessing the properties (i)--(iii) and being
different from $\cB(\cH)$, the following holds true:
$$
\cB_1(\cH)\subset \cS \subset \cB_\infty(\cH).
$$
The symmetric norm on $\cB_\infty(\cH)$ coincides with the
operator norm in $\cB(\cH)$.

Following \cite{MM99}, we recall the concept of a norm with
respect to the spectral measure of a self-adjoint operator.
\begin{definition}
\label{ENorm}
Let $Y\in \cB(\cH,\cK)$ be a bounded operator from a  Hilbert
space $\cH$ to a Hilbert space $\cK$ and let $\{E_C(\lambda)\}$ be the
spectral family of a self-adjoint {\rm(}not necessarily bounded{\rm)}
operator $C$ acting in the Hilbert space $\cK$. Introduce
\begin{equation}
\label{enorma}
\|Y\|_{E_C}=\left(\Sup\limits_{\{\delta_k\}}
\sum_k \|E_C(\delta_k)Y\|^2\right)^{1/2},
\end{equation}
where the supremum is taken over a finite {\rm(}or
countable{\rm)} system of mutually disjoint  Borel subsets
$\left\{\delta_k\right\}$, $\delta_k\cap\delta_l=\emptyset$, if
$k\neq l$.  The number $\|Y\|_{E_C}$ is called the $E_C$-norm of
the operator $Y$. For $Z\in\cB(\cK,\cH)$ the $E_C$-norm
$\|Z\|_{E_C}$ is defined as $\|Z\|_{E_C}=\|Z^*\|_{E_C}$.
\end{definition}

One easily  checks that if the norm $\|Y\|_{E_C}$ is finite one has
$$
\|Y\|\le\|Y\|_{E_C}.
$$
If, in addition, $Y$ is a Hilbert-Schmidt operator, then
\begin{equation}
\label{hsch}
\|Y\|_{E_C}\le\|Y\|_2, \quad Y\in
\cB_2(\cH, \cK),
\end{equation}
where $\|\cdot\|_2$ denotes the Hilbert-Schmidt norm in $\cB_2(\cH,\cK)$.
\begin{definition}
\label{DefSolSyl}
Let $A$ and $C$ be densely defined possibly unbounded closed operators
in the Hilbert spaces $\cH$ and $\cK$, respectively.  A bounded
operator $X\in \cB(\cH,\cK)$ is said to be a weak solution of
the Sylvester equation
\begin{equation}
\label{syl}
XA-CX=Y, \quad Y\in \cB(\cH,\cK),
\end{equation}
if
\begin{equation}
\label{sylw}
\lal XAf,g\ral-\lal Xf,C^* g\ral=\lal Yf,g\ral
\quad \text{ for all } f\in \dom(A)
\text{ and } g\in \dom(C^*).
\end{equation}
A bounded operator $X\in \cB(\cH,\cK)$ is said to be a strong
solution of the Sylvester equation \eqref{syl} if
\begin{equation}
\label{ransyl}
\ran\biggl(\reduction{X}{\dom(A)}\biggr)\subset\dom(C),
\end{equation}
and
\begin{equation}
\label{syls}
XAf-CXf=Yf \quad \text{ for all } f\in \dom(A).
\end{equation}
Finally, a bounded operator $X\in \cB(\cH,\cK)$ is said to be
an operator solution of the Sylvester equation \eqref{syl} if
$$
\ran(X)\subset \dom (C),
$$
the operator $XA$ is bounded on $\dom (XA)=\dom(A)$, and the
equality
\begin{equation}
\label{sylext}
\overline{XA}-CX=Y
\end{equation}
holds as an operator equality, where $\overline{XA}$
denotes the closure of
$XA$.
\end{definition}
Along with the Sylvester equation \eqref{syl} we also introduce
the dual equation
\begin{equation}
\label{sylZ}
ZC^*-A^*Z=Y^*,
\end{equation}
for which the notion of weak, strong, and operator solutions is
defined in a way analogous to that in Definition
\Ref{DefSolSyl}.

It is easy to see that if one of the equations \eqref{syl} or
\eqref{sylZ} has a weak solution then so does the other one.
\begin{lemma}
\label{XZsylw}
Let $A$ and $C$ be densely defined possibly unbounded closed
operators in the Hilbert spaces $\cH$ and $\cK$, respectively.
Then an operator $X\in\cB(\cH,\cK)$ is a weak solution to the
Sylvester equation \eqref{syl} if and only if the operator
$Z=-X^*$ is a weak solution to the dual Sylvester equation
\eqref{sylZ}.
\end{lemma}
\begin{proof}
According to Definition \Ref{DefSolSyl} an operator $X\in\cB(\cH,\cK)$
is a weak solution to \eqref{syl} if \eqref{sylw} holds. Meanwhile,
\eqref{sylw} implies
$$
-\lal X^* C^* g,f\ral+\lal X^*g,A^*f\ral=\lal Y^*g,f\ral
\quad \text{ for all }g\in \dom(C^*)
\text{ and } f\in \dom(A).
$$
Thus, by Definition \Ref{DefSolSyl} the operator $Z=-X^*$ is a
weak solution to the dual Sylvester equation \eqref{sylZ}. The
converse statement is proven in a similar way.
\end{proof}

The following result, first proven by M.\,Krein in 1948, gives
an ``explicit'' representation for a unique solution of the
Sylvester equation $XA-CX=Y$, provided that the spectra of
the operators $A$ and $C$ are disjoint  and one of them is a bounded
operator.  (Later, this result was independently obtained by
Y.\,Daleckii \cite{D53} and M.\,Rosenblum \cite{R56}).
\begin{lemma}
\label{Krein}
Let $A$ be a possibly unbounded densely defined closed operator
in the Hilbert space $\cH$ and $C$ a  bounded  operator in the
Hilbert space $\cK$ such that
$$
\spec(A)\cap\spec(C)=\emptyset
$$
and $Y\in  \cB(\cH,\cK)$.  Then  the Sylvester equation
\eqref{syl} has a unique operator solution
\begin{equation}
\label{ESylKrein}
X=\frac{1}{2\pi\ri}
\int_{\gamma} d\zeta\,(C-\zeta)^{-1}Y (A-\zeta)^{-1},
\end{equation}
where $\gamma$ is a union of closed contours in the complex
plane with total winding numbers $0$ around $\spec(A)$ and $1$
around $\spec(C)$ and the integral converges in the norm
operator topology.  Moreover, if $Y\in \cS$ for some symmetric
ideal $\cS\subset \cB(\cH, \cK)$ with the norm $\III\cdot \III$, then $X\in \cS$ and
$$
{\III}X{\III}\le (2\pi)^{-1}|\gamma|
\sup_{\zeta\in \gamma} \|(C-\zeta)^{-1}\|\, \|(A-\zeta)^{-1}\|\,
{\III}Y{\III},
$$
where $|\gamma|$ denotes the length of the contour $\gamma$.
\end{lemma}

If  $A$ and $C$ are  unbounded densely defined closed operators,
even with separated spectra, that is,
$
\dist\{\spec(A), \spec(C)\}>0,
$
then the Sylvester equation \eqref{syl} may not
have bounded weak solutions (see \cite{P91} for a
counterexample). Nevertheless, under some additional assumptions
equation \eqref{syl} is still weakly solvable.

The next statement is a generalization of Lemma \Ref{Krein} to
the case of unbounded operators, a result first proven by Heinz
\cite{H51}.
\begin{lemma}
\label{Kexp}
Let $A-\frac{d}{2}I$ and $-C-\frac{d}{2}I$, $d>0$, be maximal
accretive operators in  Hilbert spaces $\cH$ and $\cK$,
respectively, and $Y\in \cB(\cH, \cK)$.  Then  the Sylvester
equation \eqref{syl} has a unique weak solution
\begin{equation}
\label{XYSol1}
 X=\int_0^{+\infty}
dt\,{\rm e}^{C t}Y{\rm e}^{-A t},
\end{equation}
where the integral is understood in the weak operator topology.
Moreover, if $Y\in \cS$ for some symmetric ideal $\cS\subset
\cB(\cH, \cK)$ with the norm $\III\cdot \III$, then $X\in \cS$ and
$$
{\III}X{\III}\le \frac{1}{d}\,{\III}Y{\III}.
$$
\end{lemma}
If both $A$ and $C$ are self-adjoint operators with separated spectra one
still has a statement regarding the existence and uniqueness of
a weak solution with no additional assumptions.
\begin{theorem}
\label{Kiexp}
Let $A$ and $C$ be self-adjoint operators in Hilbert spaces
$\cH$ and $\cK$ and
\begin{equation}
\label{Tsemib}
  d=\dist\{\spec(A),\spec(C)\}>0.
\end{equation}
Then the Sylvester equation \eqref{syl} has a unique weak
solution
\begin{equation}
\label{HeinzMod}
X=\int_{-\infty}^\infty
e^{\ri tC} Y e^{-\ri tA} f_d(t)dt,
\end{equation}
where the integral is understood in the weak operator topology.
Here $f_d$ denotes any function in $L^1(\bbR)$, continuous except
at zero, such that
\begin{equation}
\label{fourie}
\int_{-\infty}^\infty e^{-\ri sx}f_d(s)ds= \frac{1}{x}
\text{ whenever } |x|\ge\frac{1}{d}\,.
\end{equation}
Moreover, if $Y\in \cS$ for some ideal $\cS\subset\cB(\cH,\cK)$
with a symmetric norm $\III\cdot \III$, then $X\in \cS$ and
\begin{equation}
\label{ozenka}
{\III} X{\III}\le\frac{c}{d}\,{\III} Y{\III},
\end{equation}
where
\begin{equation}
\label{cpi2}
c=\frac{\pi}{2}
\end{equation}
and  estimate \eqref{ozenka} with the constant $c$ given
by \eqref{cpi2} is sharp.
In particular, the estimate \eqref{ozenka}, \eqref{cpi2}
holds for any $Y\in\cB(\cH,\cK)$, that is,
\begin{equation}
\label{ozenkaR}
\|X\|\le \frac{\pi}{2d}\,\|Y\|.
\end{equation}
\end{theorem}
\begin{remark}\label{nad}
Theorem \Ref{Kiexp} with the following bounds for the best
possible constant $c$ in \eqref{ozenka}
\begin{equation}
\label{sharp}
\sqrt{\frac{3}{2}}\le c\le 2
\end{equation}
has been proven in \cite{BDMc83}.  {}From \cite{BDMc83} one can
also learn that the best possible constant in \eqref{ozenka}
admits the following estimate from above
\begin{equation}
\label{inf}
c\le\inf \big \{ \|f\|_{L^1(\bbR)}:\, f\in L^1(\bbR), \,\, \hat
f(x)=\frac{1}{x} \,\, , \,\,|x|\ge 1 \big \} ,
\end{equation}
where
$$
\hat f(x)=\int_{-\infty}^\infty e^{-\ri sx}f_d(s)ds,
\quad x\in \bbR.
$$
The fact that the infimum in \eqref{inf} equals $\pi/2$ goes
back to  B.~Sz.-Nagy and A.~Strausz {\rm(}cf. \cite{SN53}{\rm)}. The
proof of the fact that the value $c=\pi/2$ is sharp is due to R.
McEachin \cite{McE93}.
\end{remark}

The discussion of  existence of strong solutions  to the
Sylvester equation needs some technical tools from the Stieltjes
theory of integration. We briefly recall the main concepts and
results of this theory (see \cite{AG93}, \cite{AdLMSr},
\cite{MM99}, and references therein).

\begin{definition}
\label{IntDef}
Let $[a,b)\subset\bbR$, $-\infty<a<b<+\infty$.  Assume that  $C$
is a self-adjoint possibly unbounded operator in $\cK$ and
$\{E_C(\mu)\}_{\mu\in \bbR}$ is its spectral family.

The operator-valued function
$$
F:\,[a,b)\to\cB(\cK,\cH)
$$
is said to be uniformly {\rm(}resp. strongly, weakly{\rm)} integrable from the
right over the spectral measure $dE_C(\mu)$ on $[a,b)$ if the
limit
\begin{equation}
\label{defIntR}
\Int_a^b F(\mu)\,dE_C(\mu)=
\Lim_{\Max_{k=1}^n |\delta_k^{(n)}|\rightarrow0}\,\,
\sum\limits_{k=1}^n  F(\zeta_k)\, E_C(\delta_k^{(n)})
\end{equation}
exists in the uniform {\rm(}resp. strong, weak{\rm)} operator topology. Here,
$\delta_k^{(n)}=[\mu_{k-1},\mu_k)$ and
$|\delta_k^{(n)}|=\mu_k-\mu_{k-1}$,  $k=1,2,\ldots,n$, where
$a=\mu_0<\mu_1<\ldots<\mu_n=b$ is a partition of the  interval
$[a,b)$, and   $\zeta_k\in \delta_k^{(n)}$. The limit
value~{\rm(\Ref{defIntR})}, if it exists, is called the right
Stieltjes integral of the operator-valued function $F$ over the
measure $dE_C(\mu)$ on $[a,b)$.

Similarly,  the function
$$
G:[a,b)\rightarrow\cB(\cH,\cK)
$$
is said to be uniformly
{\rm(}resp. strongly, weakly{\rm)} integrable from the left over
the
measure $dE_C(\mu)$ on $[a,b)$, if there exists the
limit
\begin{equation}
\label{defIntL}
\Int_a^b dE_C(\mu)\,G(\mu)=
\Lim_{\Max_{k=1}^n |\delta_k^{(n)}|\rightarrow0}\,\,
\sum\limits_{k=1}^n E_C(\delta_k^{(n)}) G(\zeta_k)
\end{equation}
in the uniform {\rm(}resp. strong, weak{\rm)} operator topology.
The corresponding limit value {\rm(\Ref{defIntL})}, if it
exists, is called the left Stieltjes integral of the
operator-valued function $G$ over the measure $dE_C(\mu)$ on
$[a,b)$.
\end{definition}
\begin{lemma} [\rm{\cite{MM99}, Lemma 10.5}]
\label{Integr1}
An operator-valued function $F(\mu)$,
$$
F:\,[a,b)\rightarrow\cB(\cK,\cH),
$$
is integrable in the weak {\rm(}uniform{\rm)} operator topology over the
measure $dE_C(\mu)$ on $[a,b)$ from the left if and only if  the
function $[F(\mu)]^*$ is integrable in the weak {\rm(}uniform{\rm)}
operator topology over the measure $dE_C(\mu)$ on $[a,b)$ from the
right and then
\begin{equation}
\label{JJadj}
\left[\Int_a^b F(\mu)\,dE_C(\mu)\right]^*
     =\Int_a^b dE_C(\mu)\,[F(\mu)]^*.
\end{equation}
\end{lemma}
\begin{remark}
\label{RemConv}
In general, the convergence of one of the integrals
\eqref{JJadj} in the strong operator topology only implies the
convergence of the other one in the weak operator topology.
\end{remark}

Some sufficient conditions for the integrability of an
operator-valued function $F(\mu)$ over a finite interval  in the
uniform operator topology are available. For instance, we have
the following statement.
\begin{lemma}[\rm \cite{AdLMSr}, Lemma~7.2 and Remark~7.3]
\label{Integr2}
Let $\cH$ and $\cK$ be Hilbert spaces
and  let $C$ be a self-adjoint operator
in  $\cK$. Assume that the operator-valued  function
$F$, $F:[a,b)\rightarrow\cB(\cK,\cH)$,
satisfies the Lipschitz condition
\begin{equation}
\label{Lipschitz}
\|F(\mu_2)-F(\mu_1)\|\leq c\,|\mu_2-\mu_1| \quad
\text{for any}\quad \mu_1,\mu_2\in[a,b)
\end{equation}
for some constant $c>0$. Then the operator-valued function $F$
is right-integrable on $[a,b)$ with respect to the spectral
measure $dE_C(\mu)$ in the sense of the uniform operator
topology.
\end{lemma}

\medskip

The improper weak, strong, or uniform  right (left) integrals
$\Int_a^b F(\mu)\,dE_C(\mu)$ $\left (\Int_a^b
dE_C(\mu)\,G(\mu)\,\right )$ with infinite lower and/or upper
bounds ($a=-\infty$ and/or $b=+\infty$) are understood as the
limits, if they exist, of the integrals over finite intervals in
the corresponding topologies.  For example,
$$
\Int_{-\infty}^\infty dE_C(\mu)\,G(\mu)=
\Lim_{a\downarrow -\infty \,b\uparrow \infty}
\Int_{a}^{b} dE_C(\mu)\,G(\mu).
$$
We also use the notations
$$
\Int_{\spec(C)} dE_C(\mu)\,G(\mu)=
\Int_{-\infty}^{+\infty} dE_C(\mu)\,G(\mu)
$$
and
$$
\Int_{\spec(C)} F(\mu)\,dE_C(\mu)=
\Int_{-\infty}^{+\infty} F(\mu)\,dE_C(\mu).
$$
\begin{lemma}[\rm \cite{MM99}, Lemma 10.7]
\label{IntegrN}
Let an operator-valued function $F:\spec(C)\rightarrow\cB(\cH)$
be bounded
$$
\|F\|_\infty=\Sup\limits_{\mu\in\spec(C)}\|F(\mu)\|<\infty,
$$
and admit a bounded extension from $\spec(C)$ to the whole real
axis $\bbR$ which satisfies the Lipschitz
condition~{\rm(\Ref{Lipschitz})}.  If the $E_C$-norm
$\|Y\|_{E_C}$ of the operator $Y\in \cB(\cH,\cK)$ is finite,
then the integrals
$$
\Int_{\spec(C)} dE_C(\mu)\,Y\,F(\mu)
\quad
\mbox{and}
\quad
\Int_{\spec(C)} F(\mu)\,Y^*\,dE_C(\mu)
$$
exist in  the uniform operator  topology. Moreover, the following
bounds hold
\begin{equation}
\label{EstT1}
\biggl\|\Int_{\spec(C)} dE_C(\mu)\,Y\,F(\mu)\biggr\|
\leq\|Y\|_{E_C}\cdot\|F\|_\infty,
\end{equation}
\begin{equation}
\label{EstT2}
\biggl\|\Int_{\spec(C)} F(\mu)\,Y^*\,dE_C(\mu)\biggr\|
\leq\|Y\|_{E_C}\cdot\|F\|_\infty.
\end{equation}
\end{lemma}

Now we are ready to state the key result of this section:  if either $A$ or $C$ is
self-adjoint, then a strong solution to the  Sylvester equation,
if it exists, can be represented  in the form of an operator
Stieltjes integral.
\begin{theorem}
\label{SylUnb}
Let $A$ be a possibly unbounded densely defined closed operator in the Hilbert
space $\cH$ and $C$ a self-adjoint operator in the Hilbert space
$\cK$. Let $Y\in\cB(\cH,\cK)$ and suppose that $A$ and $C$ have
separated spectra, i.\,e.,
\begin{equation}\label{ACC}
\dist\{\spec(A),\spec(C)\}>0.
\end{equation}
Then the following statements are valid.
\begin{enumerate}

\item[{\rm(i)}] Assume that the Sylvester equation \eqref{syl}
has a strong solution $X\in\cB(\cH,\cK)$.  Then  $X$ is a unique
strong solution to \eqref{syl} and it can be represented in the
form of the Stieltjes integral
\begin{equation}
\label{EqX}
X=\int_{\spec(C)} E_{C}(d\mu)
Y(A-\mu)^{-1},
\end{equation}
which converges in the sense of the strong operator
topology in $\cB(\cH,\cK)$.

Conversely, if the Stieltjes integral \eqref{EqX} converges in
the strong operator topology, then $X$ given by \eqref{EqX} is a
strong solution to \eqref{syl}.

\item[{\rm(ii)}] Assume that the dual Sylvester equation
\begin{equation}
\label{sylZs}
ZC-A^*Z=Y^*
\end{equation}
has a strong solution $Z\in\cB(\cK,\cH)$.  Then  $Z$ is a unique
strong solution to~\eqref{sylZs} and it can be
represented in the form of the Stieltjes operator integral
\begin{equation}
\label{EqXZ}
Z=-\int_{\spec(C)} (A^*-\mu)^{-1} Y^* E_{C}(d\mu),
\end{equation}
which converges in the sense of the strong operator
topology in $\cB(\cK,\cH)$.

 Conversely, if the operator Stieltjes integral in \eqref{EqXZ}
converges in the strong operator topology, then $Z$ given by
\eqref{EqXZ} is a strong solution to \eqref{sylZs}.
\end{enumerate}
\end{theorem}
\begin{proof} (i)
Assume that  the Sylvester equation \eqref{syl} has a strong
solution $X\in \cB(\cH,\cK)$, that is, \eqref{ransyl} and \eqref
{syls} hold.  Let $\delta$ be a finite interval such that
$\delta\cap \text{spec }(C)\ne \emptyset$ and $\mu_\delta\in \delta\cap
\text{spec }(C)$.  Applying to both sides of \eqref{syls} the
spectral projection $E_{C}(\delta)$, a short
computation yields
\begin{equation}
\label{koki}
 E_{C}(\delta)XAf-\mu_\delta E_{C}(\delta)Xf=
 E_{C}(\delta)Yf+E_{C}(\delta)(C-\mu_\delta)Xf
\end{equation}
for any  $ f\in \dom(A)$.  Since $\mu_\delta\in \delta\cap
\text{spec }(C)$, by \eqref{ACC} one concludes that $\mu_\delta$
belongs to the resolvent set of the operator $A$. Hence,
\eqref{koki} implies
\begin{equation}\label{vagno}
 E_{C}(\delta)X=E_{C}(\delta)Y(A-\mu_\delta)^{-1}
+(C-\mu_\delta)E_{C}(\delta)X(A-\mu_\delta)^{-1}.
\end{equation}

Next, let $[a, b)$ be a finite interval and $\{\delta_k\}$ a
finite system of mutually disjoint intervals such that
$[a,b)=\cup_k \delta_k $.  For those $k$ such that $\delta_k\cap
\text{\rm spec }(C)\ne \emptyset$ pick a point
$\mu_{\delta_k}\in\delta_k\cap \text{\rm spec }(C)$.  Using
\eqref{vagno} one obtains
\begin{align}
\sum\limits_{k:\delta_k\cap\,\spec(C)\neq\emptyset}
E_{C}(\delta_k)X=&
\sum\limits_{k:\delta_k\cap\,\spec(C)\neq\emptyset}
E_{C}(\delta_k)Y
(A-\mu_{\delta_k})^{-1}
\no \\
&+\sum\limits_{k:\delta_k\cap\,\spec(C)\neq\emptyset}
(C-\mu_{\delta_k})\,E_{C}(\delta_k)
X(A-\mu_{\delta_k})^{-1}.
\label{dlinno}
\end{align}
The left hand side of \eqref{dlinno} can be computed explicitly:
\begin{equation}
\label{lili}
\sum\limits_{\delta_k\cap\,\spec(C)\neq\emptyset}
E_{C}(\delta_k)X
=E_{C}\bigl([a,b)\cap\,\spec(C)\bigr)X=
E_{C}\bigl([a,b)\bigr)X.
\end{equation}

The first term on the r.\,h.\,s. of \eqref{dlinno} is the
integral sum for the  Stieltjes integral  \eqref{EqX}.
More precisely, since $(A-\mu)^{-1}$ is  analytic in a complex
neighborhood of $[a,b]\cap\,\spec(C)$, by Lemma \Ref{Integr2}
one infers
\begin{align}
\label{Limit}
&{\nlim\limits_{\mathop{\rm max}\limits_{k}|\delta_k|\to 0}
\sum\limits_{k:\delta_k\cap\,\spec(C)\neq\emptyset}
E_{C}(\delta_k)Y
(A-\mu_{\delta_k})^{-1}}\\
\nonumber
&\qquad=\int\limits_{[a,b)\cap\,\spec(C)}
E_{C}(d\mu)Y(A-\mu)^{-1}.
\end{align}

The last term on the right hand side of \eqref{dlinno}
vanishes
\begin{equation}
\label{putu}
\nlim\limits_{\mathop{\rm max}\limits_{k}|\delta_k|\to 0}
\sum\limits_{k:\delta_k\cap\,\spec(C)\neq\emptyset}
(C-\mu_{\delta_k})\,E_{C}(\delta_k)\,
X(A-\mu_{\delta_k})^{-1}=0.
\end{equation}
This can be seen as follows.  For any $f\in\cH$ we have the
estimate
\begin{align*}
&{\biggl\|\sum\limits_{\delta_k\cap\,\spec(C)\neq\emptyset}
(C-\mu_{\delta_k})\,E_{C}(\delta_k)\,
X(A-\mu_{\delta_k})^{-1}f\biggr\|^2}\\
&\qquad=\biggl\langle
\sum\limits_{k:\delta_k\cap\,\spec(C)\neq\emptyset}
(A^*-\mu_{\delta_k})^{-1}X^*(C-\mu_{\delta_k})^2
E_{C}(\delta_k)X(A-\mu_{\delta_k})^{-1}f,f
\biggr\rangle\\
&\qquad\leq  \sum\limits_{\delta_k\cap\,\spec(C)\neq\emptyset}
|\delta_k|^2\|X\|^2\|(A-\mu_{\delta_k})^{-1}\|^2\|f\|^2\\
&\qquad\leq (b-a)\,\|X\|^2\,\|f\|^2\,\mathop{\rm max}\limits_k |\delta_k|
\Sup\limits_{\mu\in[a,b)\cap\,\spec(C)}
\|(A-\mu)^{-1}\|^2.
\end{align*}
Here we have used the estimate
$$
\|(C-\mu_{\delta_k})^2 E_{C}(\delta_k)\|=
\left\|\int_{\delta_k}(\mu-\mu_{\delta_k})^2
E_{C}(d\mu)\right\|\leq
\Sup\limits_{\mu\in\delta_k}(\mu-\mu_{\delta_k})^2\leq|\delta_k|^2.
$$
Passing to the limit  $\mathop{\rm max}\limits_k |\delta_k|\to
0$ in \eqref{dlinno}, by \eqref{lili}--\eqref{putu} one
concludes that for any finite interval $[a,b)$
\begin{equation}
\label{semf}
E_{C}\bigl([a,b)\bigr)X=
\int\limits_{[a,b)\cap\,\spec(C)}
E_{C}(d\mu)Y(A-\mu)^{-1}.
\end{equation}
Since
$$
\slim\limits_{\mbox{\scriptsize$\begin{array}{c}
a\to-\infty\\b\to+\infty
\end{array}$}}
E_{C}\bigl([a,b)\bigr)X=X,
$$
(\Ref{semf}) implies (\Ref{EqX}), which, in particular, proves
the  uniqueness of a strong solution to the Riccati equation
\eqref{syl}.

In order to prove the converse statement of (i), assume that the
Stieltjes integral on the r.\,h.\,s. part of (\Ref{semf})
converges as $a\to-\infty$ and $b\to+\infty$ in the strong
operator topology. Denote the resulting integral by $X$.  Then,
(\Ref{semf}) holds for any finite $a$ and $b$.  This implies
that for any $f\in \dom(A)$ we have
\begin{align*}
&{C E_{C}\bigl([a,b)\bigr)Xf-
E_{C}\bigl([a,b)\bigr)XA f}\\
&\qquad=\int\limits_{[a,b)\cap\,\spec(C)}
E_{C}(d\mu)Y(A-\mu)^{-1}(\mu-A)f\\
&\qquad=-\int\limits_{[a,b)\cap\,\spec(C)}
E_{C}(d\mu)Yf=-E_{C}\bigl([a,b)\bigr)Yf.
\end{align*}
Hence,
\begin{equation}\label{predric}
   C E_{C}\bigl([a,b)\bigr)Xf
=E_{C}\bigl([a,b)\bigr)XA f-
E_{C}\bigl([a,b)\bigr)Yf \text{ for any } f\in \dom(A)
\end{equation}
and $ C E_{C}\bigl([a,b)\bigr)Xf$
converges to $XA f+Y f$ as $a\to-\infty$
and $b\to+\infty$. Therefore,
$$
  \sup_{a<b} \bigl\|C E_{C}\bigl([a,b)\bigr)Xf\bigr\|^2
    =\sup_{a<b}\int\limits_{[a,b)\cap\,\spec(C)}\mu^2\,
    d\lal E_{C}Xf,Xf\ral <\infty
$$
 and, hence,
\begin{equation}
\label{XfC}
Xf\in\dom(C).
\end{equation}
Then
\eqref{predric} can be rewritten in the form
\begin{equation}\label{BB}
 E_{C} \bigl([a,b)\bigr)C Xf
=E_{C}\bigl([a,b)\bigr)XA f-
E_{C}\bigl([a,b)\bigr)Y f,\quad a<b.
\end{equation}
Combining \eqref{XfC} and \eqref{BB} proves that $X$ is a strong
solution to the Sylvester equation \eqref{syl}.


{\rm(ii)}
Assume that the dual Sylvester equation~(\Ref{sylZ}) has a
strong solution $Z\in\cB(\cK,\cH)$. As  in the proof of (i),
choose a finite interval $\delta\subset\bbR$ such that $\delta\cap
\text{spec }(C)\ne \emptyset$.  Since
$E_C(\delta)\cK\subset\dom(C)$, we have
$ZE_C(\delta)f\in\dom(A^*)$ for any $f\in\cK$ by the definition
of a strong solution.  Take a point $\mu_\delta\in \delta\cap
\text{spec }(C)$. It follows from \eqref{ACC} that
$\mu_\delta\not\in\spec(A^*)$.  As in the proof of (i), it is
easy to check the validity of the representation
\begin{equation}
\label{vagnoZ}
\begin{array}{c}
ZE_C(\delta)f=-(A^*-\mu_\delta)^{-1}Y^*E_C(\delta)f
-(A^*-\mu_\delta)^{-1}Z(C-\mu_\delta)E_C(\delta)f, \\
\end{array}
\end{equation}
which holds $\text{for all}\quad f\in\cK.$

Next, let $[a,b)$ be a finite interval and $\{\delta_k\}$ a
finite system of mutually disjoint intervals such that
$[a,b)=\cup_k \delta_k $. For those $k$ such that $\delta_k\cap
\text{\rm spec }(C)\ne \emptyset$ pick a point
$\mu_{\delta_k}\in\delta_k\cap \text{\rm spec }(C)$.
Using \eqref{vagnoZ} one then finds that
\begin{align}
\nonumber
ZE_{C}([a,b))f=&
-\sum\limits_{k:\delta_k\cap\,\spec(C)\neq\emptyset}
(A^*-\mu_{\delta_k})^{-1}Y^* E_{C}(\delta_k)f\\
\label{dlinnoZ}
&-\sum\limits_{k:\delta_k\cap\,\spec(C)\neq\emptyset}
(A-\mu_{\delta_k})^{-1}Z\,(C-\mu_{\delta_k})E_{C}(\delta_k)f.
\end{align}
The equality \eqref{putu} implies
\begin{equation}
\label{putuZ}
\nlim\limits_{\mathop{\rm max}\limits_{k}|\delta_k|\to 0}
\sum\limits_{k:\delta_k\cap\,\spec(C)\neq\emptyset}
(A^*-\mu_{\delta_k})^{-1}Z(C-\mu_{\delta_k})E_{C}(\delta_k)
=0.
\end{equation}
Thus, passing in (\Ref{dlinnoZ}) to the limit
as $\mathop{\rm max}\limits_k |\delta_k|\to 0$
one  infers that
\begin{equation}
\label{semfZ}
-\int\limits_{[a,b)\cap\,\spec(C)}
(A^*-\mu)^{-1}Y^*E_{C}(d\mu)f=ZE_{C}\bigl([a,b)\bigr)f.
\end{equation}
Since for any $f\in\cK$
$$
\mathop{\rm lim}\limits_{\mbox{\scriptsize$\begin{array}{c}
a\to-\infty\\b\to+\infty
\end{array}$}}
ZE_{C}\bigl([a,b)\bigr)f=Z,
$$
one concludes that the integral on the r.\,h.\,s. part
of~(\Ref{EqXZ}) converges as $a\to-\infty$ and $b\to+\infty$
in the strong operator topology and
(\Ref{EqXZ}) holds, which gives
 a unique strong solution to
the  dual Sylvester equation~(\Ref{sylZs}).


In order to prove   the converse statement of (ii),  assume that
there exists the strong operator limit
\begin{equation}
\label{IntZab}
Z=\slim\limits_{\mbox{\scriptsize$\begin{array}{c}
a\to-\infty\\b\to+\infty
\end{array}$}}
\int\limits_{[a,b)\cap\,\spec(C)}
(A^*-\mu)^{-1}Y^*E_{C}(d\mu) , \qquad Z\in\cB(\cK,\cH).
\end{equation}
For any finite $a$ and $b$ such that $a<b$ we have
\begin{equation}
\label{ZECab}
ZE_C([a,b)=-\int_{\spec(C)\cap[a,b)}
(A^*-\mu)^{-1} Y^* E_{C}(d\mu).
\end{equation}
By \eqref{ACC} any point $\zeta\in \spec(C)$ belongs to the
resolvent set of the operator $A$ and, hence, to the one of
$A^*$. Picking such a $\zeta$, $\zeta\in \spec(C)$, the
operator \eqref{ZECab} can be split into two parts
\begin{equation}
\label{ZEJJ}
ZE_{C}([a,b))=J_1(a,b)+J_2(a,b),
\end{equation}
where
\begin{align}
\label{J1}
J_1(a,b)&=-(A^*-\zeta)^{-1}Y^* E_{C}([a,b)),\\
\label{J2}
J_2(a,b)&=+(A^*-\zeta)^{-1}
\int\limits_{\spec(C)\cap [a,b)}(\zeta-\mu)(A^*-\mu)^{-1}Y^*
E_{C}(d\mu).
\end{align}
Using the functional calculus for the self-adjoint operator $C$
  one obtains
$$
J_2(a,b)f=
-(A^*-\zeta)^{-1}\left(\int\limits_{\,\,\spec(C)\cap [a,b)}
(A^*-\mu)^{-1}Y^* E_{C}(d\mu)\right)(C-\zeta)f,
$$
$$
\text{ for } f\in \dom (C).
$$
Thus, for $f\in \dom (C)$ one concludes that
\begin{align*}
Zf=&
\lim\limits_{\mbox{\scriptsize$\begin{array}{c}
a\to-\infty\\b\to+\infty
\end{array}$}}ZE_{C}([a,b))f
\no \\
=&\lim\limits_{\mbox{\scriptsize$\begin{array}{c}
a\to-\infty\\b\to+\infty
\end{array}$}}
J_1(a,b)f
+
\lim\limits_{\mbox{\scriptsize$\begin{array}{c}
a\to-\infty\\b\to+\infty
\end{array}$}}J_2(a,b)f
\no \\
=&-(A^*-\zeta)^{-1}Y^*f
\no \\
&-(A^*-\zeta)^{-1}\left(\int\limits_{\,\,\spec(C)}
(A^*-\mu)^{-1}Y^* E_{C}(d\mu)\right)(C-\zeta)f
\no
\end{align*}
That is,
\begin{equation}
\label{ranZA}
Zf=-(A^*-\zeta)^{-1}Y^*f+(A^*-\zeta)^{-1}Z(C-\zeta)f, \quad f\in \dom(C),
\end{equation}
since
\begin{align}
\nonumber
&{\int_{\spec(C)}
(A^*-\mu)^{-1}Y^* E_{C}(d\mu)} \\
\label{sIntDef}
&\qquad=\slim\limits_{\mbox{\scriptsize$\begin{array}{c}
a\to-\infty\\b\to+\infty
\end{array}$}}
\int\limits_{\spec(C)\cup [a,b)}
(A^*-\mu)^{-1}Y^* E_{C}(d\mu)=Z
\end{align}
by \eqref{IntZab}.
It follows  from~(\Ref{ranZA}) that
$Zf\in\dom(A^*)$ for any $f\in\dom(C)$ and, thus,
\begin{equation}
\label{ransylZ}
\ran\biggl(\reduction{Z}{\dom(C)}\biggr)\subset\dom(A^*).
\end{equation}
Applying $A^*-\zeta$ to the both sides of the resulting equality
(\Ref{ranZA}) one infers that $Z$ is a strong solution to the
dual Sylvester equation~(\Ref{sylZs}) which completes the proof.
\end{proof}

\begin{corollary}
\label{CorXZ}
Assume the hypothesis of Theorem \Ref{SylUnb}. Assume, in addition,
that the Sylvester equations \eqref{syl} has a strong solution
$X\in\cB(\cH,\cK)$. Then $Z=-X^*$ is a unique weak solution to
the dual Sylvester equation \eqref{sylZ}. Vice versa, if
$Z\in\cB(\cK,\cH)$ is a strong solution of the dual Sylvester
equation \eqref{sylZ}, then $X=-Z^*$ is a unique weak solution
to the equation \eqref{syl}.
\end{corollary}
\begin{remark}
\label{RemAsym} The proofs of  parts {\rm(i)} and {\rm(ii)}  of
Theorem \Ref{SylUnb} are slightly different in flavour owning to
the fact that the operation of taking the adjoint is not
continuous in the strong operator topology. Hence, in general, we
are not able to state that the  strong convergence of the
Stieltjes integral in \eqref{EqX} implies the strong convergence
of that in \eqref{EqXZ} and vice versa {\rm(}cf. Remark
\Ref{RemConv}{\rm)}.
\end{remark}

For the sake of completeness we also present a ``weak''  version
of Theorem \Ref{SylUnb}.
\begin{theorem}
\label{SylUnbWeak}
Assume the hypothesis of Theorem \Ref{SylUnb}. Then the following
statements are equivalent.

{\rm(i)} The Sylvester equations \eqref{syl} has a weak
solution $X\in\cB(\cH,\cK)$.

{\rm(ii)}
There exists the weak limit
\begin{equation}
\label{wlX}
X=\wlim\limits_{\mbox{\scriptsize$\begin{array}{c}
a\to-\infty\\b\to+\infty
\end{array}$}}
\int_{\spec(C)\cap[a,b)} E_{C}(d\mu)Y(A-\mu)^{-1}.
\end{equation}

{\rm(iii)} The dual Sylvester equation \eqref{sylZ} has a weak
solution $Z=-X^*\in\cB(\cK,\cH)$.

{\rm(iv)}
There exists the weak
limit
\begin{equation}
\label{wlZ}
Z=-\wlim\limits_{\mbox{\scriptsize$\begin{array}{c}
a\to-\infty\\b\to+\infty
\end{array}$}}
\int_{\spec(C)\cap[a,b)} (A^*-\mu)^{-1}Y^*E_{C}(d\mu).
\end{equation}

\end{theorem}

The statement below concerns the existence of  strong
and even operator solutions to the Sylvester equation.

\begin{lemma}
\label{SylEC}
Assume the hypothesis of Theorem \Ref{SylUnb}.
Assume, in addition, that the condition
\begin{equation}
\label{Res}
\Sup\limits_{\mu\in\,\spec(C)}
\|(A-\mu)^{-1}\|<\infty
\end{equation}
holds and the operator $Y$ has a finite $E_C$-norm, that is,
\begin{equation}
\label{VarB}
\|Y\|_{E_{C}}<\infty.
\end{equation}

Then the  Sylvester equations \eqref{syl} and \eqref{sylZ} have
unique strong and, hence, unique weak solutions
$X\in\cB(\cH,\cK)$ given by \eqref{EqX} and $Z\in\cB(\cK,\cH)$
given by \eqref{EqXZ}, respectively, and, moreover, $Z=-X^*$. In
representations \eqref{EqX} and \eqref{EqXZ} the Stieltjes
integrals exist in the sense of the uniform operator topology.

Assume,  in addition, that
\begin{equation}
\label{ResHi}
\Sup\limits_{\mu\in\,\spec(C)}
\|\mu\,(A-\mu)^{-1}\|<\infty.
\end{equation}

Then
\begin{equation}
\label{XinC}
\ran(X)\subset\dom(C),
\end{equation}
\begin{equation}
\label{ZinA}
\ran(Z)\subset\dom(A),
\end{equation}
and, thus, $X$ and $Z$ appear to be operator solutions
to \eqref{syl} and \eqref{sylZ}, respectively.
\end{lemma}
\begin{proof}
By \eqref{Res}, \eqref{VarB}, and Lemma \Ref{IntegrN}  the
operator Stieltjes integrals in \eqref{EqX} and \eqref{EqXZ} can
be understood  in the operator norm topology. Thus, $X$ given
\eqref{EqX} and $Y$ given by \eqref{EqXZ} are unique strong
solutions to the Sylvester equations \eqref{syl} and
\eqref{sylZ} by Theorem \Ref{SylUnb}.  Therefore, the  operators
$X$ and $Z$ are unique weak solutions and $Z=-X^*$ by Theorem
\Ref{SylUnbWeak}.

In order to prove \eqref{XinC} it suffices to note that under
conditions (\Ref{ResHi}) and (\Ref{VarB}) for any $f\in\cH$ and
for any $a,b\in\bbR$, $a<b$, due to \eqref{EqX} the following
estimate holds
\begin{align*}
&{\int\limits_{[a,b)\cap\,\spec(C)}\mu^2\,
d\lal E_{C}Xf,Xf\ral =
\bigl\|C E_{C}\bigl([a,b)\bigr)Xf\bigr\|^2}\\
&\qquad=\int\limits_{[a,b)\cap\,\spec(C)}
\lal Y^*E_{C}(d\mu)Y\,\mu\,(A-\mu)^{-1}f,
\,\mu\,(A-\mu)^{-1}f\ral\\
&\qquad\leq\|Y\|^2_{E_{C}}\,
\biggl(\Sup\limits_{\mu\in\,\spec(C)}
\|\mu\,(A-\mu)^{-1}\|\biggr)^2\,\|f\|^2.
\end{align*}
Thus,
$$
    \int\limits_{\spec(C)} \mu^2\,\lal E_{C}Xf,Xf\ral
    <\infty,
$$
which proves that $Xf\in\dom(C)$ and, hence, the inclusion
\eqref{XinC} is proven.

It remains to   prove the inclusion \eqref{ZinA}. Given
$\zeta\in\spec(C)$, we represent $ZE_C([a,b))$ for some finite
$a,b\in\bbR$, $a<b$, in the form \eqref{ZEJJ} where $J_1(a,b)$
and $J_2(a,b)$ are just the same ones  as in (\Ref{J1})
and (\Ref{J2}), respectively.  Under condition \eqref{ResHi}, by
Theorem \Ref{IntegrN} one concludes that the operator Stieltjes
integral in~(\Ref{J2}) converges as $a\to-\infty$ and
$b\to+\infty$ in the uniform operator topology to some operator
$M\in\cB(\cK,\cH)$.  Then, from \eqref{ZEJJ} one learns that for
any $f\in\cK$
\begin{align*}
Zf&=
\lim\limits_{\mbox{\scriptsize$\begin{array}{c}
a\to-\infty\\b\to+\infty
\end{array}$}}ZE_{C}([a,b))f
\no \\
&=\lim\limits_{\mbox{\scriptsize$\begin{array}{c}
a\to-\infty\\b\to+\infty
\end{array}$}}
J_1(a,b)f
+
\lim\limits_{\mbox{\scriptsize$\begin{array}{c}
a\to-\infty\\b\to+\infty
\end{array}$}}J_2(a,b)f
\no \\
&=-(A^*-\zeta)^{-1}Y^*f+(A^*-\zeta)^{-1}Mf
\no
\end{align*}
and, thus, $Zf\in\dom(A)$ which proves \eqref{ZinA}.

The proof is complete.
\end{proof}

\begin{remark}
If the operator $A$ is self-adjoint, then the strong solution
of the Sylvester equation, if it  exists, can be represented in
the form of the repeated Stieltjes integral
\begin{equation}
\label{rep}
X=\int_{\spec(C)}dE_C(\mu)Y\int_{\spec(A)}\frac{dE_A(\lambda)}
{\lambda -\mu}.
\end{equation}
If, in addition, $Y$ is a Hilbert-Schmidt operator, then the repeated
integral \eqref{rep} can also be represented in the form of the
double Stieltjes integral
\begin{equation}
\label{doubl}
X=\int\int_{\spec(C)\times\spec(A)}\frac{dE_C(\mu)YdE_A(\lambda)}
{\lambda -\mu},
\end{equation}
where the integral \eqref{doubl} can be understood as the
$\cB_2$-norm limit of the integral sums of the Lebesgue type. It
is also worth to mention that by a theorem by Birman and
Solomjak \cite{BS67} under condition \eqref{ACC} we have the
estimate
\begin{equation}
\label{hssh}
\|X\|_2\le \frac{1}{d}\, \|Y\|_2,
\end{equation}
where  $ d=\dist\{\spec(A),\spec(C)\}.
$
Moreover, the estimate \eqref{hssh} is sharp in the class of
Hilbert-Schmidt operators.
\end{remark}
\begin{remark}
If $Y$ is a Hilbert-Schmidt operator, inequality \eqref{hssh} is
a considerable improvement of the more general estimate
\eqref{ozenka}, the latter being sharp only in the class of all
symmetric normed ideals. We also remark that if $A$ is
self-adjoint and \eqref{VarB} holds, then \eqref{EqX} implies
the estimate
\begin{equation}
\label{XYEC}
\|X\|_{E_C}\leq\frac{1}{d}\,\|Y\|_{E_C}.
\end{equation}
\end{remark}

\section{Riccati equation}
\label{SecRic}

The goal of this section is to develop an approach for solving
the operator Riccati equations based on an applications of
Banach's Fixed Point Principle for transformations of operator
spaces. Putting aside the discussion of the purely geometric
approach suggested  and developed by  Davis and Kahan
\cite{DK69}, \cite{DK70} and by Adams \cite{A83} as well as the
one based on the factorization technique for operator
holomorphic functions by Markus and  Matsaev \cite{MM89},
\cite{MarkusMatsaev75} (see also \cite{MenShk}, \cite{MM99},
\cite{VirozubMatsaev}, and \cite{LR}) we concentrate ourselves
on applications of a purely analytic approach based on the
representation theorems of Section \Ref{SecSylv}.

\begin{definition}
\label{DefSolRic}
Assume that  $A$ and $C$ are  possibly unbounded densely defined
closed operators in the Hilbert spaces $\cH$ and $\cK$,
respectively, while $B\in\cB(\cK,\cH)$ and $D\in\cB(\cH,\cK)$.

A bounded operator $Q\in\cB(\cH,\cK)$ is said to be a weak
solution of the Riccati equation
\begin{equation}
\label{RicABCD}
QA-CQ+QBQ=D
\end{equation}
if
$$
%
\begin{array}{c}
\lal QAf,g\ral-\lal Qf,C^*g\ral+\lal QBQf,g\ral=\lal Df,g\ral\\[0.5em]
 \text{ for all } f\in \dom (A)
\text{ and } d\in \dom(C^*).
\end{array}
$$

A bounded operator $Q\in\cB(\cH,\cK)$ is said to be a strong
solution of the Riccati equation \eqref{RicABCD} if
\begin{equation}
\label{ranric}
\ran\biggl(\reduction{Q}{\dom(A)}\biggr)\subset\dom(C),
\end{equation}
and
\begin{equation}
\label{rics}
QAf-CQf+QBQf=Df  \text{ for all } f\in \dom(A).
\end{equation}

Finally, a bounded operator $Q\in \cB(\cH,\cK)$ is said to be
an operator solution of the Riccati equation \eqref{RicABCD} if
$$
\ran(Q)\subset \dom (C),
$$
the operator $QA$ is bounded on $\dom(QA)=\dom(A)$
and the equality
\begin{equation}
\label{Ricext}
\overline{QA}-CQ+QBQ=D
\end{equation}
holds as an operator equality, where $\overline{QA}$ denotes the closure of
$QA$.
\end{definition}

Along with the Riccati equation \eqref{RicABCD} we also introduce
the dual equation
\begin{equation}
\label{RicK}
KC^* - A^*K + KB^*K=D^*,
\end{equation}
for which the notion of weak, strong, and operator solutions
is defined in a way analogous to that in Definition
\Ref{DefSolRic}.
\begin{example}
\label{Examp0}
{\rm(The Friedrichs model~\cite{Fried})}.  Given a nonempty open
Borel set $\Delta\subset \bbR$, let  $\cH=\bbC$ and
$\cK=L^2(\Delta)$. Let  $A=0$ in $\cH$ and let $C$ be the
multiplication operator in $\cK$,
$$
(Cf)(\mu)= \mu f(\mu)
$$
on
$$
\dom(C)=\{f\in L^2(\Delta):\, \int_\Delta\,d\mu\,
(1+\mu^2)|f(\mu)|^2<\infty,
$$
$B\in\cB(\cK,\cH)$, and, finally, $D=B^*\in\cB(\cH,\cK)$.
\end{example}

By Riesz representation theorem
$$
Bf=\lal f,b\ral=\int_{\Delta}\, d\mu \, f(\mu)\overline{b(\mu)},
\quad f\in \cK,
$$
for some essentially bounded function $b\in \cK=L^2(\Delta)$ and
hence
$$
(D\zeta) (\mu)=\overline{ b(\mu)}\zeta, \quad \zeta\in \bbC,
$$
since $D=B^*$.

Under the assumptions of this example a bounded operator
$Q\in \cB(\cH, \cK)$ is a weak  solution to  the Riccati
equation \eqref{RicABCD} if and only if $Q$ has the form
\begin{equation}
\label{form}
(Q\zeta)(\mu)=q(\mu)\zeta, \quad\zeta\in\bbC,
\end{equation}
where $q$ is an essentially bounded function, and
\begin{equation}
\label{Ricq}
-(\mu q)(\mu)+\lal q,b\ral q(\mu)=b(\mu)
\quad \text{ for a.\,e. } \mu\in\Delta.
\end{equation}
Moreover, any weak solution $Q$ appears to be a strong solution,
that is, any essentially bounded function $q$ satisfying
\eqref{Ricq} belongs to $\dom (C)$.

Solving \eqref{Ricq} with respect to $q$ one concludes that the
Riccati equation \eqref{Ricq} has a weak/strong solution if and
only if
\begin{equation}
\label{L2}
\text{there exists a $w\in\bbR$ such that  }
\frac{b(\cdot)}{\cdot -w}\in L^2(\Delta)
\end{equation}
and
\begin{equation}
\label{wbasic}
w+\int_{\Delta}d\mu\, \frac{\,\,|b(\mu)|^2}
{\mu-w}=0.
\end{equation}
If conditions \eqref{L2} and \eqref{wbasic} hold for some
$w\in\bbR$, then  the solution $Q$ has the form \eqref{form},
where
\begin{equation}
\label{qbasic}
q(\mu)=\frac{b(\mu)}{w-\mu},\quad\mu\in\bbR,
\end{equation}
and
$$
w=\lal q,b\ral.
$$
\medskip

The next assertion is a direct corollary of Lemma
\Ref{XZsylw}.
\begin{lemma}
\label{QKricw}
Let $A$ and $C$ be densely defined possibly unbounded
closed operators in the Hilbert spaces $\cH$ and $\cK$,
respectively, $B\in\cB(\cK,\cH)$, and $D\in\cB(\cH,\cK)$.
Then $Q\in\cB(\cH,\cK)$ is a weak solution to the Riccati equation
\eqref{RicABCD} if and only if $K=-Q^*$ is a weak
solution to the dual Riccati equation \eqref{RicK}.
\end{lemma}

Throughout the remaining part of the section we assume the
following hypothesis.
\begin{hypothesis}
\label{hhqq}
Assume that $\cH$ and $\cK$ are Hilbert spaces, $A$ and $C$ are
possibly unbounded self-adjoint operators on domains $\dom(A)$
in $\cH$ and $\dom(C)$ in $\cK$, respectively.  Also assume that
$B\in\cB(\cK,\cH)$ and $D\in\cB(\cH,\cK)$.
\end{hypothesis}

The representation theorems of  Sec.\,\Ref{SecSylv} for solutions
of the Sylvester equation are a source for iteration schemes
which allow one to prove solvability of Riccati equations by
using fixed point theorems.  Here we present two of such
schemes for the search for strong or weak solutions to the
Riccati equation.
\begin{theorem}
\label{basics}
Assume Hypothesis \Ref{hhqq}.  Then the following statements
hold true.

\begin{enumerate}

\item[{\rm(i)}] Assume, in addition to Hypothesis \Ref{hhqq},
that
$$
  \dist\{\spec(A),\spec(C)\}>0.
$$
Then $Q\in\cB(\cH,\cK)$ is a weak solution to the
Riccati equation \eqref{RicABCD} if and only if it is a solution
to the equation
\begin{equation}
\label{FHeinz0}
Q=\int_{-\infty}^{\infty} e^{\ri tC}\,(D-QBQ)\,e^{-\ri tA} f_d(t) dt,
\end{equation}
where $f_d$ is a summable function satisfying  \eqref{fourie}
and the integral in \eqref{FHeinz0} exists in the sense of
the weak operator topology in $\cB(\cH,\cK)$.

\item[{\rm(ii)}] Assume, in addition to Hypothesis \Ref{hhqq}, that
\begin{equation}
\label{ABQCsep}
\dist\{\spec(A+BQ),\spec(C)\}>0.
\end{equation}
Then $Q\in\cB(\cH,\cK)$ is a strong {\rm(}weak{\rm)} solution to
the Riccati equation \eqref{RicABCD} if and only if $Q$ is a
solution of the equation
\begin{equation}
\label{EqQ}
Q=\int_{\spec(C)} E_{C}(d\mu)
D\,(A+BQ-\mu)^{-1},
\end{equation}
where  the operator Stieltjes integral exists in the sense of
the strong {\rm(}weak{\rm)} operator topology in $\cB(\cH,\cK)$.

\item[{\rm(iii)}]
Assume, in addition to Hypothesis \Ref{hhqq}, that
$K\in\cB(\cK,\cH)$ and
\begin{equation}
\label{ABKCsep}
\dist\{\spec(A-KB^*),\spec(C)\}>0.
\end{equation}
Then the operator $K$ is a strong {\rm(}weak{\rm)} solution to
the dual Riccati equation \eqref{RicK} if and only if $K$
satisfies the equation
\begin{equation}
\label{BasicK}
K=-\int_{\spec(C)} (A-KB^*-\mu)^{-1}D^*E_C(d\mu),
\end{equation}
where the operator Stieltjes integral exists in the sense of the
strong {\rm(}weak{\rm)} operator topology.

\end{enumerate}
\end{theorem}
\begin{proof} (i) The operator $Q$ is a weak solution to
\eqref{RicABCD} if and only if $Q$ is a weak solution to the
equation
$$
QA-CQ=Y,
$$
where
$$
Y=D-QBQ.
 $$
Applying Theorem \Ref{Kiexp} completes the proof of (i).

(ii)  The operator $Q$ is a strong solution to
\eqref{RicABCD} if and only if $Q$ is a strong solution to the
equation
$$
Q\tA-CQ=D,
$$
where
$$
\tA=A+BQ.
$$
Applying Theorem \Ref{SylUnb} (i) completes the proof of (ii).

(iii)
 The operator $K$ is a strong solution to
\eqref{RicK} if and only if $K$ is a strong solution to the
equation
$$
KC-\widehat AK=D^*,
$$
where
$$
\widehat A=A-KB^*.
$$
Applying Theorem \Ref{SylUnb} (ii) completes the proof of (iii).

The proof is complete.
\end{proof}

The following statement is a direct consequence of Lemma
\Ref{SylEC}.
\begin{theorem}
\label{RicEC}
Assume Hypothesis \Ref{hhqq} and
let  $D$ have  a finite norm with respect to
the spectral measure of the operator $C$, that is,
\begin{equation}
\label{VarDC}
\|D\|_{E_{C}}<\infty.
\end{equation}
Assume, in addition, that an operator $Q\in\cB(\cH,\cK)$
is a weak solution of the Riccati equation \eqref{RicABCD} such that
\begin{equation}
\label{ABQCsep1}
\dist\{\spec(A+BQ),\spec(C)\}>0,
\end{equation}
and  that the condition
\begin{equation}
\label{Ress}
\Sup\limits_{\mu\in\,\spec(C)}
\|(A+BQ-\mu)^{-1}\|<\infty
\end{equation}
holds.

Then  $Q$ is a strong solution to
\eqref{RicABCD} and  the operator $K=-Q^*$
is a strong solution to the dual Riccati equation
\eqref{RicK}.

 The strong solutions $Q$ and $K$ admit the representations
 \begin{equation}
\label{EqQ1}
Q=\int_{\spec(C)} E_{C}(d\mu)
D\,(A+BQ-\mu)^{-1},
\end{equation}
\begin{equation}
\label{BasicK1}
K=-\int_{\spec(C)} (A-KB^*-\mu)^{-1}D^*E_C(d\mu),
\end{equation}
where the operator Stieltjes integrals
 exist in the sense of the uniform operator
topology.
 Hence, the operators $Q$ and $K$ have finite
$E_C$--norm and the following bound holds true
\begin{equation}
\label{QEC0}
\|K\|_{E_C}=\|Q\|_{E_C}\leq \|D\|_{E_C}
\Sup_{\mu\in\spec(C)}\|(A+BQ-\mu)^{-1}\|.
\end{equation}
If, in this case, instead of \eqref{Ress} the following condition
holds
\begin{equation}
\label{ResHiQ}
\Sup\limits_{\mu\in\,\spec(C)}
\|\mu\,(A+BQ-\mu)^{-1}\|<\infty,
\end{equation}
then
$$
\ran(Q)\subset\dom(C)
$$
 and
 $$
 \ran(K)\subset\dom(A)
 $$
 and,
hence,
 the strong solutions $Q$ and $K$ appear to be
operator solutions to the Riccati
equations \eqref{RicABCD} and \eqref{RicK}, respectively.
\end{theorem}

In the case where the spectra of the operators $A$ and $C$ are
separated, under additional ``smallness" assumptions upon the
operators $B$ and $D$ we are able to prove the existence of
fixed points for mappings given by \eqref{FHeinz0} and
\eqref{EqQ}.
\begin{theorem}
\label{QsolvN}
Assume Hypothesis \Ref{hhqq}
and suppose that
$$
   B\neq0.
$$
Also assume that
\begin{equation}
\label{tsemib}
  d=\dist\{\spec(A),\spec(C)\}>0.
\end{equation}
Then:

\begin{enumerate}

\item[{\rm(i)}] If the inequality holds
\begin{equation}
\label{BDpi2}
\sqrt{\|B\|\,\|D\|}<\frac{d}{\pi},
\end{equation}
then the Riccati equation \eqref{RicABCD} has a unique weak
solution in the ball
$$
\left\{Q\in\cB(\cH,\cK):\, \|Q\|<\frac{d}{\pi\|B\|}\right\}.
$$
The weak solution $Q$ satisfies the estimate
\begin{equation}
\label{EstL}
\|Q\|\leq \frac{1}{\|B\|}\left(\frac{d}{\pi}-\sqrt{\frac{d^2}{\pi^2}-
\|B\|\,\|D\|}\right).
\end{equation}
In particular, if
\begin{equation}
\label{BDpiSum}
\|B\|+\|D\|<\frac{2}{\,\pi\,}\,d,
\end{equation}
then the weak solution $Q$ is a strict contraction,
that is,
$$
\|Q\|<1.
$$

\item[{\rm(ii)}] If the operator $D$ has a finite $E_C$--norm and the
inequality
\begin{equation}
\label{BDest}
\sqrt{\|B\|\,\|D\|_{E_C}}<\frac{d}{2}
\end{equation}
holds, then the Riccati equation \eqref{RicABCD} has a unique strong solution
in the ball
\begin{equation}
\label{QEst2}
\left\{Q\in\cB(\cH,\cK):\,
\|Q\|<\|B\|^{-1}\left(d-\sqrt{\|B\|\,\|D\|_{E_C}}\right)\right\}.
\end{equation}
The strong solution $Q$ has a finite $E_C$--norm and one has
the estimate
\begin{equation}
\label{QEst1}
\|Q\|_{E_C}\leq
 \frac{1}{\|B\|}\,
\left(\frac{d}{2}-\sqrt{\frac{d^2}{4}-\|B\|\,\|D\|_{E_C}}\right).
\end{equation}
In particular, if
\begin{equation}
\label{BDestSum}
\|B\|+\|D\|_{E_C}<d,
\end{equation}
then the strong solution  $Q$ is a strict contraction in
both the uniform operator and $E_C$-norm topologies, that is,
$$
\|Q\|\le \|Q\|_{E_C}<1.
$$

\end{enumerate}

\end{theorem}
\begin{proof}
The proof is based on an application of Banach's Fixed Point
Theorem.

{\rm(i)}   Let $f\in L^1(\bbR)$ be
a continuous function on $\bbR$ except at zero
such that
$$
\widehat f(s)=\int_\bbR e^{-\ri st}f(t)dt=\frac{1}{s}
\,\,\, \text{ whenever } |s|\ge 1.
$$
Introducing the function
$$
f_d(t)=f(dt), \quad t\in \bbR,
$$
by Theorem \Ref{basics} {\rm(i)} any fixed point of the map $F(Q)$
given by
\begin{equation}
\label{FHeinz}
F(Q)=\int_{-\infty}^\infty
e^{\ri tC} (D-QBQ) e^{-\ri tA} f_d(t)dt, \quad Q\in\cB(\cH,\cK),
\end{equation}
where the improper Riemann integral is understood in the weak
sense, is a weak solution to the Riccati equation
\eqref{RicABCD}.  Taking into account  that
$$
\|f_d\|_{L^1(\bbR)}=\frac{\|f\|_{L^1(\bbR)}}{d},
$$
 from \eqref{FHeinz} one concludes that
\begin{equation}
\label{EstF1}
\|F(Q)\|\leq \frac{
\|f\|_{L^1(\bbR)}
}{d}(\|D\|+\|B\|\,\|Q\|^2),
 \quad Q\in\cB(\cH,\cK)
\end{equation}
and
\begin{equation}
\label{EstF2}
\|F(Q_1)-F(Q_2)\|\leq \frac{\|f\|_{L^1(\bbR)}}{d}
\|B\|\,(\|Q_1\|+\|Q_2\|)\, (\|Q_1-Q_2\|),
\end{equation}
$$
\quad Q_1,Q_2\in\cB(\cH,\cK).
$$

Clearly, $F$ maps the ball
$\cO_r=\{Q\in\cB(\cH,\cK):\, \|Q\|\leq r\}$
into itself whenever
$$
\frac{\|f\|_{L^1(\bbR)}}{d}(\|D\|+\|B\|\,r^2)\le r
$$
and $F$ is a strict contraction of the ball  $\cO_r$ whenever
$$
\frac{2 \|f\|_{L^1(\bbR)}\|B\|}{d}\,r<1.
$$

Since the extremal problem for the Fourier transform,  which is
to find the infimum of $\|f\|_{L^1}$ over all functions $f\in
L^1(\bbR)$ such that $\widehat f(s)=1/s$ for $|s|\ge 1$, has the
solution (cf. Remark \Ref{nad})
$$
\inf \,\{ \|f\|_{L^1(\bbR)}\, :\, f\in L^1(\bbR), \widehat
f(s)=1/s \,\,\,\text{ whenever } |s|\ge 1 \}=\frac{\pi}{2},
$$
one concludes that $F$ maps the ball $\cO_r=\{Q\in\cB(\cH,\cK):\,
\|Q\|\leq r\}$ into itself whenever
\begin{equation}
\label{eq1}
\frac{\pi}{2d}(\|D\|+\|B\|\,r^2)\le r
\end{equation}
and $F$ is a strict contraction of $\cO_r$
whenever
\begin{equation}
\label{eq2}
\frac{\pi \|B\|}{d}r<1.
\end{equation}

Solving inequalities \eqref{eq1} and \eqref{eq2} one concludes that
if the radius  $r$ of the ball $\cO_r$ is within the bounds
\begin{equation}
\label{IneqR}
\frac{d}{\pi \|B\|}-\sqrt{\frac{d^2}{
\pi^2\|B\|^2}- \frac{\|D\|}{\|B\|}} \leq r <
\frac{d}{\pi\|B\|},
\end{equation}
then $F$ is a strictly contractive mapping of the ball $\cO_r$ into
itself. Applying Banach's Fixed Point Theorem  proves
assertion (i).

{\rm(ii)} Given $r\in (0,d\,\|B\|^{-1})$, under Hypothesis
\eqref{tsemib} we have the identity
\begin{align}
\label{ResIdentity}
(A+BQ-\mu)^{-1}=&
\left(I+(A-\mu)^{-1}BQ\right)^{-1}(A-\mu)^{-1},\\
\nonumber
&\mu\in \spec(C),
\,\,Q\in \cO_r,
\end{align}
which implies  the estimate
\begin{align}
\nonumber
\sup_{\mu\in \spec(C)}\|(A+BQ-\mu)^{-1}\| \leq& \displaystyle
\sup_{\mu\in \spec(C)}\frac{1}{1-\|(A-\mu)^{-1}\|\,\|B\|\,\|Q\|}\,
\|(A-\mu)^{-1}\|\\
\label{ResExt}
\leq& \displaystyle\frac{1}{1-\frac{\|B\| r}{d}}\,
\frac{1}{d}=\displaystyle\frac{1}{d-\|B\| r},
\end{align}
whenever
$
 Q\in \cO_r.
$

Since \eqref{ResExt} holds and the operator $D$ has a finite
$E_C$-norm,
the mapping
$$
F(Q)=\int_{\spec(C)} E_{C}(d\mu)D\,(A+BQ-\mu)^{-1},
$$
where the integral is understood in the strong sense, is well
defined on the domain
$$
\dom(F)=\cO_r.
$$
Since for $Q\in \cO_r$ one clearly has the estimate
$$
\dist\{\spec(A+BQ),\spec(C)\}\geq d-\|B\|r >0,
$$
any fixed point of the map $F$ is a strong solution to the
Riccati equation \eqref{RicABCD} by Theorem \Ref{basics}
{\rm(ii)}.

Using \eqref{ResExt}
we have the following  two estimates
\begin{align}
\|F(Q)\|\leq\|F(Q)\|_{E_C}&\leq \|D\|_{E_C}
\Sup_{\mu\in\spec(C)}\|(A+BQ-\mu)^{-1}\|
\no\\
&\le\frac{\|D\|_{E_C}}{d-\|B\|r}, \quad Q\in \cO_r,
\label{QQFC}
\end{align}
and
\begin{align}
\no
&\|F(Q_1)-F(Q_2)\|\\
&\qquad\leq\|F(Q_1)-F(Q_2)\|_{E_C}
\no \\
&\quad\qquad=\left\|\,\int_{\spec(C)} E_C(d\mu)D\,
 (A+BQ_1-\mu)^{-1}\,B(Q_2-Q_1)\,(A+BQ_2-\mu)^{-1}\,
\right\|_{E_C}
\no \\
&\quad\qquad\leq\frac{\|D\|_{E_C}}{(d-\|B\|r)^2}
\|Q_2-Q_1\|,
\label{EEST}\quad Q_1,Q_2\in\cO_r.
\end{align}
Clearly, by \eqref{QQFC} $F$ maps the ball $\cO_r$ into itself
whenever
\begin{equation}
\label{Est1}
\frac{\|D\|_{E_C}}{d-\|B\|r}\leq r
\end{equation}
and by \eqref{EEST} $F$ is a strict contraction on
$\cO_r$ whenever
\begin{equation}
\label{Est2}
\frac{\|D\|_{E_C}}{(d-\|B\|r)^2}<1\,.
\end{equation}
Solving inequalities \eqref{Est1} and \eqref{Est2}
simultaneously, one concludes that if the radius of the ball
$\cO_r$ is within the bounds
\begin{equation}
\label{Br}
\frac{1}{\|B\|}\left(\frac{d}{2}-\sqrt{\frac{d^2}{4}-\|B\|\,\|D\|_{E_C}}
\right)
\leq  r < \frac{1}{\|B\|}
\left(d-\sqrt{\|B\|\,\|D\|_{E_C}}
\right),
\end{equation}
then $F$ is a strictly contracting mapping of the ball $\cO_r$
into itself.  Applying Banach's Fixed Point Theorem we infer
that equation (\Ref{EqQ}) has a unique solution in any ball
$\cO_r$ whenever $r$ satisfies \eqref{Br}.  Therefore, the fixed
point does not depend upon the radii satisfying~(\Ref{Br}) and
hence it belongs to the smallest of these balls. This
observation proves the estimate
\begin{equation}
\label{QnEst}
\|Q\|\leq\frac{1}{\|B\|}\left(
\frac{d}{2}-\sqrt{\frac{d^2}{4}-\|B\|\,\|D\|_{E_C}}\right).
\end{equation}
Finally, using
(\Ref{QQFC}),  for the fixed point  $Q$ one
obtains the estimate
\begin{equation}
\label{QEC1}
\|Q\|_{E_C}=\|F(Q)\|_{E_C}\leq
\frac{\|D\|_{E_C}}{d-\|B\|\,\|Q\|}.
\end{equation}
Then~(\Ref{QnEst}) yields
$$
\|Q\|_{E_C}\leq\frac{\|D\|_{E_C}}
{\frac{d}{2}+\sqrt{\frac{d^2}{4}-\|B\|\,\|D\|_{E_C}}}=
\frac{1}{\|B\|}
\left({\frac{d}{2}-\sqrt{\frac{d^2}{4}-\|B\|\,\|D\|_{E_C}}}\right),
$$
which
completes the proof.
\end{proof}
\begin{remark}
Part {\rm(ii)} of the theorem extends results obtained in
\cite{MotSPbWorkshop}, \cite{MotRem}, and \cite{MM99}.  In case
where the self-adjoint operator $C$ is bounded, $D$ is a
Hilbert-Schmidt operator, B is bounded, and $A$ is possibly unbounded densely defined
closed non-self-adjoint operator, the solvability of the
equation \eqref{EqQ1} under condition \eqref{tsemib} has
recently been studied in \cite{AdLT}.
\end{remark}
\begin{remark}
\label{conti}
Under the hypotheses \eqref{tsemib}  and \eqref{BDpi2} or \eqref{BDest}
the fixed point   $Q$ depends continuously {\rm(}in the operator
norm{\rm)} upon the operators $B$ and $D$, which follows from a
result {\rm(}see, e.\,g., \cite{KA84} {\rm Ch. XVI, Theorem 3)}
concerning the continuity of the mapping in Banach's Fixed
Point Theorem with respect to a parameter.
\end{remark}
\begin{remark}
\label{RemNomit}
In general,   hypothesis \eqref{tsemib} in Theorem \Ref{QsolvN}
can not be omitted.  In order to see this assume that
$\Delta=\bbR$ in Example \Ref{Examp0} and, thus,  \eqref{tsemib}
does not hold.  Assume, in addition, that the function
$b(\cdot)$ in this example is a strictly positive continuous function.
Then the necessary condition \eqref{L2} for the solvability of
the Riccati \eqref{Ricq} is violated.
\end{remark}

In order to complete the discussion of the results of Theorem
\Ref{QsolvN} we need the following illustrative statement based
on Example \Ref{Examp0}.
\begin{lemma}
\label{show}
Assume the hypothesis of Example \Ref{Examp0} for
$$
\Delta=(-\infty,-d)\cup(d,+\infty)
$$ and  some $d>0$.

If $b\in L^2(\Delta)$ and
\begin{equation}
\label{sqrt}
\|b\|\le \sqrt{2}d,
\end{equation}
then the Riccati equation \eqref{Ricq} has a weak/strong
solution. Moreover, the constant $\sqrt{2}$ in \eqref{sqrt} is
sharp.
\end{lemma}
\begin{proof}
Under the hypothesis \eqref{sqrt} we have the inequalities
$$
 \int_{-\infty}^{-d}d\mu\,
 \frac{\,\,|b(\mu)|^2}
{\,d-\mu\,}< \frac{\|b\|^2}{2d}\le d \qquad \text{and} \qquad
\int_{d}^{+\infty}d\mu\, \frac{\,\,|b(\mu)|^2}
{\,\mu+d\,}< \frac{\|b\|^2}{2d}\le d\, .
$$
The Herglotz function
$$
f(w)=w+\int_{\Delta}d\mu\, \frac{\,\,|b(\mu)|^2}
{\,\mu-w\,}
$$
is  a  strictly increasing continuous function on $(-d,d)$ and
\begin{align}
f(-d-0)=\lim_{\varepsilon \downarrow 0}
f(-d-\varepsilon)&\le
-d+\int_d^{+\infty}d\mu\, \frac{\,\,|b(\mu)|^2}
{\,\mu+d\,}<0 ,
\no \\
f(d+0)=\lim_{\varepsilon \downarrow 0}f(d+\varepsilon) &\ge
d+\int_{-\infty}^{-d}d\mu\, \frac{\,\,|b(\mu)|^2}
{\,\mu-d\,}>0,
\no
\end{align}
not withstanding the possibility for  the one-sided limits
$f(-d-0)$ and  $f(d+0)$ to turn into $-\infty$ and $+\infty$,
respectively. Therefore, the equation
$$
f(w)=0
$$
has a unique root $w_0\in (-d,d)$, the function
$$
q(\mu)=\frac{b(\mu)}{\mu-w_0}, \quad \mu \in \Delta,
$$
is an element of $L^2(\Delta)$,  and, hence, the Riccati
equation \eqref{Ricq} has a weak/strong solution, since the
existence criterion \eqref{L2}, \eqref{wbasic} is satisfied.

In order to prove that the constant $\sqrt{2}$ in the upper
bound \eqref{sqrt} is sharp, it suffices to show that for
any $c>1$ there exists a function $b\in L^2(\Delta) $ such that
$$
\|b\|=\sqrt{2}c\,d
$$
and the Riccati equation \eqref{Ricq} has no solutions $q\in
L^2(\Delta)$.

Let $\omega\in L^1(\bbR_+)$ be a positive continuous function on
$[0, \infty)$ such that
$$
\int_0^\infty \omega(t)dt =1.
$$
Given $\varepsilon >0$,  introduce the functions
$$
\omega_\varepsilon(t)=\varepsilon^{-1}\omega(t/\varepsilon),
\quad t \ge 0,
$$
and
\begin{equation}
\label{bspec}
\varphi_\varepsilon (\mu)=\left\{\begin{array}{ll}
\displaystyle\,\arctan(d+\mu)
\omega_\varepsilon^{1/2}(d-\mu),& \mu\le-d, \\
\displaystyle
\omega_\varepsilon^{1/2}(\mu-d),& \mu\ge d  .  \\
\end{array}\right .
\end{equation}
One infers
$$
\lim_{\varepsilon \downarrow 0} \| \varphi_\varepsilon\|^2=1
$$
and
\begin{equation}
\label{bol}
\lim_{\varepsilon \downarrow 0}
\left(\int_{-\infty}^{-d}
\frac{\,\,|\varphi_{\varepsilon}( \mu)|^2}{\mu+d}d\mu+
\int_{d}^{+\infty} \frac{\,\,|\varphi_{\varepsilon} (\mu)|^2}{\mu+d}
d\mu\right)
=\frac{1}{2d} .
\end{equation}

Hence for any $c>1$, one can find an $\varepsilon_0 >0$ such that the
following inequality holds
\begin{equation}
\label{bod}
\|\varphi_{\varepsilon_0} \|^{-2}
\left (
\int_{-\infty}^{-d} \frac{\,\,|\varphi_{\varepsilon_0} ( \mu)|^2}{\mu+d}d\mu+
\int_{d}^{+\infty} \frac{\,\,|\varphi_{\varepsilon_0} (\mu)|^2}{\mu+d}d\mu
\right )
>\frac{1}{2dc^2}\,\,.
\end{equation}

Introducing
$$
b(\mu)=\sqrt{2}cd\,
\frac{\varphi_{\varepsilon_0}(\mu)}{
\|\varphi_{\varepsilon_0}\|}, \quad \mu\in\Delta= (-\infty, -d]\cup [d, \infty) ,
$$
one obviously concludes that
$$
 \|b\|=\sqrt{2}cd.
 $$
Meanwhile,  \eqref{bod} implies the estimate
$$
\int_{-\infty}^{-d} \frac{\,|b ( \mu)|^2}{\mu+d}d\mu +
\int_{d}^{+\infty} \frac{|b(\mu)|^2}{\mu+d}d\mu
>d.
$$
Therefore, the Herglotz function $f(w)$ given by
$$
f(w)=w+\int_{-\infty}^{-d} \frac{\,\,|b(\mu)|^2}{\mu-w}d\mu+
\int_{d}^{+\infty} \frac{\,\,|b(\mu)|^2}{\mu-w}d\mu
$$
does not vanish on $[-d,d)$ (note that $f(w)\to+\infty$ as
$w\uparrow d$) and hence \eqref{wbasic} is violated for all
$w\in [-d,d)$.  Since $b(\cdot)$ is a continuous  function and it
does not vanish on $(-\infty, d) \cup [d, \infty)$, the
condition \eqref{L2} is violated for all
$w\in (-\infty, d)\cup [d, \infty)$.
Hence,  the Riccati equation \eqref{Ricq} has
no weak/strong solutions in this case since the existence
criterion \eqref{L2}, \eqref{wbasic} is violated.
\end{proof}
\begin{remark}
\label{bound}
The result of Lemma \Ref{show} combined with that of Theorem
\Ref{QsolvN} shows the following.

\begin{enumerate}

\item[{\rm(i)}] There is a constant $c>0$ such
that the conditions \eqref{tsemib} and
\begin{equation}
\label{BDpi3}
\|B\|\,<c\dist\{\spec(A),\spec(C)\},
\end{equation}
imply the existence of a weak solution to the Riccati equation
$$
QA-CQ+QBQ=B^*\,.
$$

\item[{\rm(ii)}] In general, the ``smallness" requirement on $B$
\eqref{BDpi3} can not be omitted {\rm(}cf. \eqref{BDpi2} and
\eqref{BDest}{\rm)}.

\item[{\rm(iii)}] The sharp value of the constant $c$ in
\eqref{BDpi3} is within the bounds
$$
\frac{1}{\pi}\le c \le \sqrt{2}.
$$

\end{enumerate}
\end{remark}

\section{The spectral shift function}
\label{SecSSF}

The main purpose of this section is to recall the concept of the
spectral shift function associated with a pair of self-adjoint
operators and to extend this concept to the case of pairs of closed
operators that are similar to self-adjoint operators.

The spectral shift function $\xi(\lambda, H,A)$ for a pair of
self-adjoint operators $(H,A)$ in a Hilbert space $\cH$ is
usually associated with the Lifshits-Krein trace formula
\begin{equation}
\label{5.1}
\tr(\varphi(H)-\varphi(A) )=\int_{\bbR}
d\lambda\, \varphi'(\lambda)\,\xi(\lambda, H,A).
\end{equation}
The trace formula \eqref{5.1} holds for a wide class of functions
$\varphi:\,\bbR\to\bbC$, including $C_0^\infty(\bbR)$,
provided that the self-adjoint operators $H$ and $A$ are
resolvent comparable, that is,
\begin{equation}
\label{nucl}
 (H-z)^{-1}-(A-z)^{-1}\in \cB_1(\cH),\quad \Img(z)\ne 0.
\end{equation}
If \eqref{nucl} holds, then we will also say that $H$ is a relatively trace
class perturbation of $A$.

The trace formula \eqref{5.1} determines the spectral shift
function up to an arbitrary complex constant. This constant
may, however, be chosen in such a way that makes the spectral
shift function to be real-valued.

In case of trace class perturbations, i.\,e., if
$$
\overline{H-A}\in \cB_1(\cH),
$$
the additional requirement that
$$
\xi(\cdot, H,A)\in L^1(\bbR)
$$
determines the spectral shift function uniquely. Being chosen in
this way, the spectral shift function $\xi(\lambda, H,A)$ can be
computed by Krein's formula via the perturbation determinant
\begin{align}
\label{krdet}
\xi(\lambda,H,A)=&\pi^{-1} \lim_{\varepsilon \downarrow 0}
\text{\rm arg} \det
\left((H-\lambda-\ri\varepsilon)
(A-\lambda-\ri\varepsilon)^{-1}\right)\\
\nonumber
&\quad \text{ for a.\,e. } \lambda \in \bbR.
\end{align}

In the case of relatively trace class perturbations
(\Ref{nucl}),  there is in general no natural way to choose
the above constant uniquely. Moreover, the
requirement of continuity of the spectral shift function
$\xi(\lambda,H,A)$ in an appropriate function space topology
with respect to small deviations of the end points $H$ and $A$,
leads to the  conclusion:  the spectral shift function
$\xi(\lambda, H,A)$ can not be introduced uniquely as a function
of the pair $(H,A)$.  It should be considered to be either  a
function of a continuous path (in an appropriate operator
topology) connecting the end points $H$ and $A$ from the same
connected component, or to be a path independent but a
multi-valued function of the spectral parameter.

However, one can introduce  the spectral shift function uniquely
modulo $\bbZ$ in such a way that for any pairs $(H, A)$,
$(H,{\tH})$, and $({\tH},A)$ of self-adjoint
operators $A$, $H$, and ${\tH}$ in $\cH$,
satisfying \eqref{nucl}, the following chain rule holds
(see \cite{Ya92})
\begin{equation}
\label{chain}
\xi(\lambda, H,A)=
 \xi(\lambda, H,{\tH})+\xi(\lambda,{\tH},A)\quad
(\modulo\bbZ)
\quad \text{for a.\,e. }\lambda \in \bbR.
\end{equation}

The extension of a concept of a spectral shift function to the
case of operators similar to self-adjoint needs additional
considerations.

We start with a definition of a zero trace commutator class.

\begin{definition}
\label{DefTC}
Let $\cA(\cH)$ be the set of all bounded
operators $V\in\cB(\cH)$ possessing the property:
\begin{equation}
\label{VR0}
\tr \big (VR-R\,V \big )=0
\end{equation}
whenever $R\in\cB(\cH)$ and
\begin{equation}
\label{VR}
VR-R\,V\in\cB_1(\cH).
\end{equation}
The set $\cA(\cH)$
is called the {\em zero trace commutator class}.
\end{definition}

In the case of an infinite-dimensional Hilbert space $\cH$
\eqref{VR} does not imply \eqref{VR0} in general.  For example,
let $P$ be a one-dimensional orthogonal projection.  Then there
is a partial isometry $S$ such that $SS^*=I$ and $S^*S=I-P.$
Taking $R=S^*$ and $V=S$ one obtains $VR-RV=P\in \cB_1(\cH),$
but $\tr\big(VR-RV\big)=1$, and, thus, \eqref{VR0} fails despite
\eqref{VR} holds true. Therefore, the zero trace commutator
class $\cA(\cH)$ is a proper subset of $\cB(\cH)$ if the Hilbert
space $\cH$ if infinite-dimensional.

\begin{lemma}
\label{dixi}
Assume that $R,V\in\cB(\cH)$ and
at least one of the following conditions holds:
\begin{enumerate}
\item[(i)] $V\in \cB_1(\cH)$;
\item[(ii)] $VR$ and $RV$ are trace class operators;
\item[(iii)] $V$ is a normal operator and $R\in \cB_2(\cH)$;
\item[(iv)] $V$ is a self-adjoint operator
            and $R\in \cB_\infty(\cH)$;
\item[(v)] $V$ is a self-adjoint operator having no absolutely
           continuous spectral subspaces;
\item[(vi)] $V$ is a normal operator with
             purely point spectrum;
\item[(vii)] $R\in \cB_p(\cH)$ and $V\in \cB_q(\cH)$ with
            $\frac{1}{p}+\frac{1}{q}=1$.
\end{enumerate}
Then $\tr(VR-RV )=0$ whenever  $VR-RV$ is a trace class operator.
\end{lemma}
\begin{remark}
The part {\rm(i)} is obvious.  The part {\rm(ii)}  follows from
Lidskii's theorem.  The statement {\rm(iii)} is due to G.~Weiss
\cite{W78}.  Assertion  {\rm(iv)} has been proven by J.~Helton
and R.~Howe \cite{HH75}. The part {\rm(v)} immediately follows
from a result by R.~W.~Carey and J.~D.~Pincus \cite{CP76} which
states that any self-adjoint operator having no absolutely
continuous spectral subspace is the sum of an operator with
purely point spectrum and a trace class one with arbitrary small
trace norm\footnote{We are indebted to Vadim Kostrykin who has
attracted our attention to this fact.}.  The results {\rm(vi)}
and {\rm(vii)} have recently been proven by V.~Lauric and
C.~M.~Pearcy \cite{LP00}.
\end{remark}

Lemma \Ref{dixi} shows that the zero trace commutator class
$\cA(\cH)$ is a rather rich set. In particular, $\cA(\cH)$
contains all the trace class operators, that is,
$$
\cB_1(\cH)\subset \cA(\cH).
$$
More generally, any operator of the form
$$
\widehat V=V+T, \quad V \in \cA(\cH), \quad T\in \cB_1(\cH),
$$
is an element of $\cA(\cH)$.  The class $\cA(\cH)$ also contains
all normal bounded operators $V$ with purely point spectrum and
all self-adjoint bounded operators having no absolutely
continuous spectrum and, therefore, in this case,
 if $V\in \cA(\cH)$ and $V$ has a
bounded inverse, then $V^{-1}\in \cA(\cH)$ as well.

\begin{definition}
\label{defadmis}
Let $H$ be a  possibly unbounded densely defined closed operator
in $\cH$ on $\dom(H)$ with $\spec(H)\subset\bbR$. The operator
$H$ is said to be admissible if there exists a self-adjoint
operator $\widehat H$  such that
\begin{enumerate}
\item[(i)]
$H$ is similar to  $\widehat{H}$, i.\,e.,
$$
H=V^{-1} \widehat{H} V \text{ on } \dom(H)=V^{-1}(\dom(\widehat{H}))
$$
for some $V\in \cA(\cH)$ such that  $V^{-1}\in \cB(\cH)$;
\item[(ii)] $H$ and  $\widehat H$ are resolvent comparable, i.\,e.,
\begin{equation}
\label{simA}
(H-z)^{-1}-(\widehat H-z)^{-1}\in \cB_1(\cH),\quad \Img(z)\ne0.
\end{equation}
\end{enumerate}
We will call the operator $\widehat H$ a self-adjoint
representative of the admissible operator $H$.
\end{definition}
Clearly, any self-adjoint operator is admissible. Moreover,
 an admissible operator may have different self-adjoint
representatives.
\begin{lemma}\label{ochevidno}
Let $H$ be an admissible operator and $\widehat H$ any
self-adjoint representative of $H$. Then
\begin{equation}\label{nuli}
\tr\big ((H-z)^{-1}-(\widehat H-z)^{-1}\big )=0, \quad \Img(z)\ne0.
\end{equation}
\end{lemma}
\begin{proof}
By the definition of an admissible operator,  the difference of
the resolvents of $H$ and $\widehat H$ is a trace class operator
and the following representation holds for some $V\in \cA(\cH)$
such that $V^{-1}\in \cB(\cH)$
\begin{align}
&(H-z)^{-1}-(\widehat H-z)^{-1}=
V^{-1}(\widehat H-z)^{-1}V-(\widehat H-z)^{-1}
\no \\
&
=[V^{-1}(\widehat H-z)^{-1}]V-V[V^{-1}(\widehat H-z)^{-1}]\in \cB_1(\cH),\quad \Img(z)\ne0,
\label{resolv}
\end{align}
which implies \eqref{nuli}, since $V\in \cA(\cH)$.
\end{proof}

\begin{corollary}\label{amis}
Let $H $  be an  admissible operator in $\cH$
and $\widehat H_1$ and $\widehat H_2$ its self-adjoint representatives
from Definition
\Ref{defadmis}. Then  $\widehat H_1$ and $\widehat H_2$
 are resolvent comparable and
\begin{equation}\label{ksidlia}
\xi(\lambda;\widehat H_1 ,\widehat H_2 )=0\quad (\modulo\bbZ)
\quad \text{for a.\,e. }\lambda \in \bbR,
\end{equation}
where   $\xi(\lambda;\widehat H_1 ,\widehat H_2)$ is the spectral
shift function associated with the pair $(\widehat H_1,\widehat H_2)$ of
self-adjoint operators.
\end{corollary}

Now we are ready to extend the concept of the spectral shift
function to the case of pairs of admissible operators.
\begin{definition}
\label{defSSF1}
Let $(H,A)$ be a pair of resolvent comparable
 admissible operators in $\cH$ and $(\widehat H, \widehat A)$ a pair of their
 self-adjoint representatives from Definition \Ref{defadmis}.
Define the spectral shift function $\xi(\lambda; H,A)$
associated with the pair $(H,A)$ by
$$
\xi(\lambda; H,A)=\xi(\lambda;\widehat{H}, \widehat A)\quad (\modulo\bbZ)
\quad \text{for a.\,e. }\lambda \in \bbR,
$$
where   $\xi(\lambda;\widehat{H},\widehat A)$ is the spectral
shift function associated with the pair $(\widehat{H},\widehat A)$ of
self-adjoint operators.
\end{definition}

The  result of Corollary \Ref{amis} combined with the chain rule
\eqref{chain} for the pairs of self-adjoint operators shows that
the spectral shift function associated with a  pair $(H,A)$ of
resolvent comparable admissible operators  is well-defined
modulo $\bbZ$, that is, it is independent of the choice of the
self-adjoint representatives $\widehat H$ and $\widehat A$ for
the  operators $H$ and $A$, respectively. In particular, we have
arrived at the following  result.
\begin{lemma}\label{stability}
Assume that   $\widehat H$ and $\widehat A$ are
self-adjoint operators and $V, V^{-1}\in \cB(\cH)$ such that
\begin{align}
\label{aklass}
&V\in \cA(\cH),\\
\label{abtr}
&(\widehat H-z)^{-1}-(\widehat A-z)^{-1}\in \cB_1(\cH),
\quad \Img(z)\ne 0,\\
\label{vabtr}
&(V^{-1}\widehat HV-z)^{-1}-(\widehat A-z)^{-1}\in \cB_1(\cH),
\quad \Img(z)\ne 0,
\end{align}
then {\em the stability property} holds
\begin{equation}
\label{stab}
\xi(\lambda;V^{-1}\widehat HV,\widehat A)=
\xi(\lambda;\widehat H,\widehat A)
\quad (\modulo\bbZ)
\quad \text{for a.\,e. }\lambda \in \bbR.
\end{equation}
\end{lemma}

The next example shows that the requirements \eqref{abtr} and
\eqref{vabtr} by themselves do not imply \eqref{stab}, if
condition \eqref{aklass} is violated.
\begin{example}
\label{example}
Let $H$ be the closure of the operator $H_0=-\frac{d^2}{dx^2}$
on  $L^2(\bbR)$ initially defined on the domain
$\dom(H_0)=C_0^\infty(\bbR)$ and $\widehat{H}$  the operator
which acts in $L^2((-\infty, 0))\oplus L^2((0, \infty))$ and
corresponds to the Dirichlet boundary condition at zero. The
difference $(H-z)^{-1}-(\widehat{H}-z)^{-1}$, $\Img(z)\ne 0$, is
rank one and, therefore, $H$ is a relatively trace class
perturbation of $\widehat{H}$.  The operators $H$ and $\widehat
H$ are obviously unitary equivalent and, therefore, there exists
a unitary operator $V$ such that  $\widehat{H}=V^{*}HV$. The
spectral shift function associated with the pair
$(\widehat{H},H)$ is known \cite{GS96} to be a half on the
essential spectrum and zero otherwise,
$$
\xi(\lambda,\widehat{H}, H)
=\frac{1}{\,2\,}\chi_{[0, \infty)}(\lambda)\quad
(\modulo\bbZ)\quad
\text{ for a.\,e. } \lambda\in \bbR\, ,
$$
where $\chi_{\Delta}(\lambda)$ denotes the characteristic
function of the Borel set $\Delta$.  Therefore,
\begin{equation}
\label{notst}
0=\xi(\lambda, H,H)\ne \xi(\lambda, V^*HV, H)=
\frac{1}{\,2\,}\chi_{[0, \infty)}(\lambda)
\end{equation}
on a set of positive Lebesgue measure. Representation
\eqref{notst} shows that  the stability property \eqref{stab}
for the spectral shift function does not hold in this case.
\end{example}

The concept of a spectral shift function associated with a pair
of admissible operators turns out to be  rather useful in the
context of not only additive but also multiplicative theory of
perturbations. The following  theorem illustrates such an
application to the multiplicative theory of perturbations in
case where the spectral shift function can be computed via
the perturbation determinant.  The corresponding representation
appears to be an immediate analog of Krein's formula
\eqref{krdet} in the self-adjoint case. The precise statement is
as follows.
\begin{theorem}
\label{main2}
Let $A$  be  a possibly unbounded self-adjoint operator in $\cH$
with domain $\dom(A)$, $B=B^*$  a trace class self-adjoint
operator,
\begin{equation}
\label{BBB}
B\in \cB_1(\cH),
\end{equation}
and  $V$  a bounded operator with a bounded inverse such that
\begin{equation}
\label{IV}
I-V\in \cB_1(\cH).
\end{equation}
Assume, in addition,  that
\begin{equation}
\label{ranV}
\ran(I-V)\subset \dom(A),
\end{equation}
the domain $\dom (A)$ is $V$-invariant
\begin{equation}
\label{Vinv}
V\dom(A)=\dom(A),
\end{equation}
and the commutator $AV-VA$, initially defined on $\dom(A)$, is a
closable operator and its closure is a trace class operator,
that is,
\begin{equation}
\label{Vcom}
\overline{AV-VA}\in \cB_1(\cH).
\end{equation}
Then for  the operator $H$ defined by
\begin{equation}
\label{defH}
H=V^{-1}(A+B)V\quad\text{on}\quad \dom(H)=\dom(A)
\end{equation}
the following holds true.
\begin{enumerate}

\item[{\rm(i)}]  The operator  $H$
is admissible. Moreover, the spectral shift function
$\xi(\lambda;H,A)$ is well defined and
\begin{equation}\label{HVA}
\xi(\lambda;H, A)=\xi(\lambda; A+B,A)
 \quad (\modulo\bbZ), \quad
\text{ for a.\,e. } \lambda \in \bbR.
\end{equation}

\item[{\rm(ii)}]
\begin{equation}
\label{HAB}
(H-z)(A+B-z)^{-1}-I\in \cB_1(\cH), \quad \Img(z)\ne 0,
\end{equation}
and, hence, the perturbation determinant
$$
D_{H/(A+B)}(z)=\det\left((H-z)(A+B-z)^{-1}\right),\quad \Img(z)\ne 0,
$$
is well defined
and, moreover,
$$
D_{H/(A+B)}(z)=1.
$$

\item[{\rm(iii)}] The perturbation determinant $D_{H/A}(z)$ is
well defined
$$
D_{H/A}(z)=\det((H-z)(A-z)^{-1}),\quad \Img(z)\ne 0,
$$
and the  spectral shift function for the admissible pair $(H,A)$
can be computed via the perturbation determinant as follows
\begin{align}
\label{det}
\xi(\lambda;H, A)=&\pi^{-1}
\lim_{\varepsilon \downarrow  0}
\mathop{\text{\rm arg}}
D_{H/A}(\lambda+\ri\varepsilon)\quad (\modulo\bbZ)\\
\nonumber
&\text{ for a.\,e. } \lambda \in \bbR.
\end{align}
\end{enumerate}
\end{theorem}
\begin{proof}
(i) Hypothesis \eqref{IV} implies that

a) $V\in \cA(\cH)$

\noindent and

b) the operators $H$ and $A+B$ are resolvent comparable.

Thus, $H$ is an admissible operator. By \eqref{BBB} the
operator $A+B$ is a trace class perturbation of $A$ and hence
$H$ and $A$ are resolvent comparable. Therefore, \eqref{HVA}
holds by the definition of the spectral shift function for a
pair of resolvent comparable admissible operators, which proves
(i).

(ii) We start with the representation
$$
(A+B-z)V(A+B-z)^{-1}=I+W(z),
\quad \Img(z)\ne 0,
$$
where
\begin{equation}
\label{tut}
W(z)=(A+B-z)(V-I)(A+B-z)^{-1},\quad \Img(z)\ne 0,
\end{equation}
makes sense by \eqref{ranV}.  By \eqref{tut}
\begin{align*}
W(z)=&V-I+(AV-VA)(A+B-z)^{-1}\\
   &+(BV-VB)(A+B-z)^{-1}\,,
\quad\Img(z)\ne 0,
\end{align*}
which proves that
\begin{equation}
\label{tut1}
W(z)\in \cB_1(\cH),
\quad\Img(z)\ne 0,
\end{equation}
by \eqref{BBB},  \eqref{IV} and \eqref{Vcom}.  Therefore, the
Fredholm determinant of the operator $(A+B-z)V(A+B-z)^{-1}$ is
well defined and
\begin{align}
\no
&\det\big ((A+B-z)V(A+B-z)^{-1}\big)\\
\label{repr}
&\qquad=\det \big
(I+(A+B-z)(V-I)(A+B-z)^{-1}\big ),
\quad \Img(z)\ne 0.
\end{align}
Since \eqref{ranV} holds, the operator $(A+B-z)(V-I)$ is well
defined on the whole Hilbert space $\cH$ as a closed operator
being the product of two closed operators. Hence, $(A+B-z)(V-I)$
is bounded by the Closed Graph Theorem.  In particular,  the
following representation holds
\begin{equation}
\label{repr1}
(A+B-z)^{-1}[(A+B-z)(V-I)]=V-I.
\end{equation}
Using \eqref{repr},  \eqref{repr1}, and
the fact that
$$
\det(I+ST)=\det(I+TS), \quad ST, \, TS \in \cB_1(\cH),
$$
one proves
\begin{equation}
\label{repr4}
\det((A+B-z)V(A+B-z)^{-1})=\det(V),
\quad\Img(z)\ne 0.
\end{equation}
Further, using  definition \eqref{defH} of $H$
one computes
\begin{equation}
\label{det1}
(H-z)(A+B-z)^{-1}=V^{-1}(A+B-z)V(A+B-z)^{-1}),\quad\Img(z)\ne 0,
\end{equation}
which  by \eqref{tut} and \eqref{tut1} proves \eqref{HAB}.
Moreover, \eqref{det1} and \eqref{repr4} yield
\begin{align*}
\det((H-z)(A+B-z)^{-1})=&\det (V^{-1})\det(A+B-z)V(A+B-z)^{-1})\\
=&\det (V^{-1})\det (V)=1,\quad\Img(z)\ne 0,
\end{align*}
which completes the proof of (ii).

(iii)
One infers
$$
(H-z)(A-z)^{-1}=(H-z)(A+B-z)^{-1}(A+B-z)(A-z)^{-1},
\quad\Img(z)\ne 0.
$$
Hence
$$
(H-z)(A-z)^{-1}-I\in \cB_1(\cH),\quad\Img(z)\ne 0,
$$
by \eqref{HAB} and the fact that
$$
(A+B-z)(A-z)^{-1}-I\in \cB_1(\cH),\quad\Img(z)\ne 0,
$$
since $B\in  \cB_1(\cH)$, which proves that the perturbation
determinant
$D_{H/A}(z)$  is well defined.
Moreover,
\begin{equation}
\label{detdet}
D_{H/A}(z)=D_{H/(A+B)}(z)D_{(A+B)/A}(z)=D_{(A+B)/A}(z),
\quad\Img(z)\ne 0.
\end{equation}
By Krein's formula \eqref{krdet} we have
$$
\xi(\lambda;A+B, A)=\pi^{-1}
\lim_{\varepsilon \downarrow   0}
\mathop{\text{\rm arg}} D_{(A+B)/A}(\lambda+\ri\varepsilon)
\quad (\modulo\bbZ),
$$
and hence
\eqref{det} holds by \eqref{detdet}.
\end{proof}
\begin{remark}
The idea of introducing the spectral shift function associated
with a pair of operators  similar to self-adjoint operators via
the perturbation determinant (in the framework of the trace
class perturbations theory) goes back to V. Adamjan and H.
Langer \cite{AdL00}.  The proof of Theorem \Ref{main2} contains
some fragments of their original reasoning.
\end{remark}

\section{Graph subspaces and block diagonalization of operator
matrices}
\label{SecGPP}

In this section we collect some results related to existence of
invariant graph subspaces of a linear operator and to the
closely related problem of block diagonalization of the operator
in terms of such subspaces.

First, we recall the  definition of a graph subspace.
\begin{definition}
\label{DefGraf}
Let $\cN$  be a closed subspace of a Hilbert space $\cH$
and $Q\in \cB({\cN},{\cN}^\perp)$. The set
$$
 \cG({\cN},Q)=\{x\in{\cH}:\,
P_{{\cN}^\perp}\,x=QP_{\cN}\,x\}
$$
is called the graph subspace of ${\cH}$ associated with the pair
$({\cN},Q)$, where
$P_{\cN}$ and $P_{{\cN}^\perp}$  denote the
orthogonal projections  onto ${\cN}$ and
${\cN}^\perp$, respectively.
\end{definition}
It is easy to check that
\begin{equation}
\label{GGs}
\cG({\cN},Q)^\perp=
\cG({\cN}^\perp,-Q^*).
\end{equation}

{}From the analytic point of view, the search for
invariant/reducing graph subspaces for a linear self-adjoint
operator in $\cH$  is equivalent to the problem of solving the
operator Riccati equations studied in details in Section
\Ref{SecRic}.

We adopt the following hypothesis in the sequel.
\begin{hypothesis}
\label{diag}
Assume that the Hilbert space $\cH$ is decomposed into the
orthogonal sum of two orthogonal subspaces
\begin{equation}
\label{decom}
\cH=\cH_0\oplus\cH_1,
\end{equation}
the self-adjoint operator  $\bH$ reads with respect to the
decomposition \eqref{decom} as a $2\times2$ operator block
matrix
\begin{equation}
\label{twochannel}
{\bH}=\left(\begin{array}{lr}
  A_0      &   B_{01}  \\
  B_{10}   &   A_{1}
\end{array}\right),
\end{equation}
where $A_i$, $i=0,1$, are self-adjoint operators in $\cH_i$ with
domains $\dom(A_i)$ while $B_{ij}\in\cB(\cH_j,\cH_i)$, $j=1-i$,
are bounded operators and $B_{10}=B_{01}^*$.  Thus,
\begin{align}
\bH &=\bA+\bB,
\label{bAbB} \\
\dom(\bH) &=\dom(\bA),
\end{align}
where $\bA$  is the diagonal self-adjoint operator,
\begin{eqnarray}
\label{bA}
{\bA}&=&\left(\begin{array}{lr}
  A_0        & 0    \\
  0   &   A_{1}
\end{array}\right), \\
\nonumber
\dom(\bA)&=&\dom(A_0)\oplus\dom(A_1),
\end{eqnarray}
and the operator $\bB=\bB^*$ is an off-diagonal bounded operator
\begin{equation}
\label{bB}
{\bB}=\left(\begin{array}{lr}
  0        & B_{01}    \\
  B_{10}   &     0
\end{array}\right).
\end{equation}
\end{hypothesis}

We start with a criterion of existence of the
invariant graph subspaces $\cG(\cH_i,Q_{ji})$
($Q_{ji}\in\cB(\cH_i,\cH_j)$), $i=0,1,$  $j=1-i$, associated
with the $2\times2$ block decomposition \eqref{twochannel} of a
self-adjoint operator $\bH$.
\begin{lemma}
\label{princip}
Assume Hypothesis \Ref{diag}. The graph subspace
$\cG_i=\cG(\cH_i,Q_{ji})$ for some  $Q_{ji}\in\cB(\cH_i,\cH_j)$,
$i=0,1,$ $j=1-i$, is a reducing subspace for the operator $\bH$
if and only if the operator Riccati equation
\begin{equation}
\label{riceq}
\bQ\bA-\bA\bQ+\bQ\bB\bQ=\bB,
\end{equation}
has a strong solution $\bQ$ which reads with
respect to the decomposition~\eqref{decom} as
\begin{equation}
\label{QQ}
\bQ=\left(\begin{array}{cc}
0      &   Q_{01}\\
Q_{10} &   0
\end{array}\right)
\end{equation}
with
\begin{equation}
\label{QQs}
Q_{01}=-Q_{10}^*.
\end{equation}
\end{lemma}
\begin{proof}
If $\bQ$ given by \eqref{QQ}, \eqref{QQs} is a strong solution of
\eqref{riceq}, this means that
\begin{equation}
\label{ranq}
\ran\left(\reduction{\bQ}{\dom(\bA)}\right)\subset\dom(\bA),
\end{equation}
and
\begin{equation}
\label{riceqst}
\bQ\bA f-\bA\bQ f+\bQ\bB\bQ f=\bB f
\quad \text{ for any } f\in \dom(\bA).
\end{equation}

Under  hypotheses \eqref{QQ}, \eqref{QQs}, and \eqref{ranq} we
have the inclusions
\begin{equation}
\label{dodo}
\ran\left(\reduction{Q_{ji}}{\dom(A_i)}\right)\subset\dom(A_j),
\quad i=0,1,\quad j=1-i.
\end{equation}
Moreover, the Riccati equation  \eqref{riceqst} splits into a
pair of the equations
\begin{equation}
\label{figas}
Q_{ji}A_if-A_j Q_{ji} f+Q_{ji}
B_{ij}Q_{ji} f = B_{ji}f \quad \text{ for all } f\in \dom(A_i),
\end{equation}
$$
\quad i=0,1,\quad j=1-i.
$$
Rewriting these equations in the form
\begin{equation}
\label{invg}
Q_{ji}(A_i+B_{ij}Q_{ji})f=(B_{ji}+A_j Q_{ji}) f
\quad \text{ for all } f\in \dom(A_i)
\end{equation}
one immediately observes that \eqref{invg} combined with
\eqref{dodo} is equivalent to invariance of the subspaces
$\cG_i=\cG(\cH_i,Q_{ji})$, $i=0,1,$ $ j=1-i,$ for the operator $\bH$. In turn,
\eqref{QQs} implies the
invariance of the subspace $\cG_i^\perp=\cG(\cH_j,-Q_{ji}^*),$
$i=0,1,$ $ j=1-i,$
for $\bH$, which proves the lemma.
\end{proof}
\begin{remark}
\label{notalways}
Example \Ref{Examp0}  shows that, in general, the Riccati
equations \eqref{figas} are not always solvable and, thus, the
invariant graph subspaces may not  always exist either.
\end{remark}

If the operator block matrix $\bH$ has reducing graph subspaces,
then the block diagonalization problem can be solved explicitly.
\begin{theorem}
\label{thHi2}
Assume Hypothesis \Ref{diag}. Assume, in addition, that the
graph subspaces $\cG_i=\cG(\cH_i,Q_{ji})$  for some
$Q_{ji}\in\cB(\cH_i,\cH_j)$, $i=0,1,$ $ j=1-i,$  satisfying
\eqref{QQs} are reducing subspaces for the operator $\bH$. Then:
\begin{enumerate}
\item[(i)] The operator $\bV=\bI+\bQ$ with $\bQ$ given by \eqref{QQ},
           \eqref{QQs} has a bounded inverse.
\item[(ii)] The operator $\bV^{-1}\bH\bV$ is block diagonal with
           respect to the decomposition \eqref{decom}.  That is,
\begin{equation}
\label{btH}
\bV^{-1}\bH\bV=\left(\begin{array}{cc}
A_0+B_{01}Q_{10}     &   0\\
0 &  A_1+B_{10}Q_{01}
\end{array}\right)
\end{equation}
where
\begin{eqnarray}
\label{Hi}
\dom(A_i+B_{ij}Q_{ji})=\dom(A_i), \quad i=0,1, \quad j=1-i.
\end{eqnarray}
\item[(iii)] The operator $\bU^*\bH\bU$, where $\bU$ is the
     unitary operator from the polar decomposition
     $\bV=\bU|\bV|$, is  block-diagonal with respect to
     the decomposition \eqref{decom}.  That is,
\begin{equation}\label{UH}
\bU^*\bH\bU=\left(\begin{array}{cc}
H_0      &   0\\
0 &   H_1
\end{array}\right)
\end{equation}
with
\begin{equation}\label{Hip}
H_i=(I_{\cH_i}+Q^*_{ji}Q_{ji})^{1/2}(A_i+B_{ij}Q_{ji})
(I_{\cH_i}+Q^*_{ji}Q_{ji})^{-1/2},
\end{equation}
$$
 i=0,1,\quad
j=1-i,
$$
\begin{equation}\label{Hi'dom}
\dom(H_i)=(I_{\cH_i}+Q^*_{ji}Q_{ji})^{1/2}(\dom(A_i)),
\end{equation}
where $I_{\cH_i}$ stands for the identity operator in $\cH_i$.
\end{enumerate}
\end{theorem}
\begin{proof}

(i) By \eqref{QQs} $\bQ^*=-\bQ$ and, thus, the spectrum of $\bQ$
is a subset of the imaginary axis. This means that zero does not
belong to the spectrum of $\bV=\bI+\bQ$ and, hence, $\bV$ has a
bounded inverse.

(ii) Since by (i) $\bV$ has a bounded inverse, \eqref{btH} is
equivalent to the representation
$$
\bH\bV=\bV\left(\begin{array}{cc}
A_0+B_{01}Q_{10}      &   0\\
0 &   A_1+B_{10}Q{01}
\end{array}\right),
$$
which, in turn, taking into account \eqref{Hi},
 is equivalent to the Riccati equation
\eqref{riceq}. Then, applying Lemma \Ref{princip}, the validity
of \eqref{btH}--\eqref{Hi} is equivalent to the fact that the
graph subspaces $\cG_i=\cG(\cH_i,Q_{ji})$, $i=0,1,$ $j=1-i,$
are reducing subspaces.

(iii) Taking into account \eqref{QQs}, by inspection one gets
\begin{equation}
\label{bV2}
\bV\bV^*=\bV^*\bV=
\left(\begin{array}{cc}
I_0+Q_{01}Q_{01}^* & 0\\
0& I_1+Q_{10}Q_{10}^*
\end{array}\right).
\end{equation}
Since $\bV=\bU|\bV|$ and $ |\bV|=(\bV\bV^*)^{1/2}$, the validity
of \eqref{UH}--\eqref{Hi'dom} is an immediate consequence of
\eqref{btH}--\eqref{Hi}.
\end{proof}

\section{Invariant graph subspaces and
splitting of the spectral shift function}
\label{SecKey}

It is convenient to
study spectral properties of the perturbed block operator matrix
$\bH$ not only in terms of the perturbation $\bB=\bH-\bA$ in
itself, but also in terms of the angular operator $\bQ$
associated with the reducing graph subspaces, provided that they
exists.  The next (conditional) result throws  light upon the
{\it quantitative} aspects of the perturbation theory for block operator
matrices in this context.

\begin{theorem}
\label{th42}
Assume Hypothesis \Ref{diag} and let  the Riccati equation
\eqref{riceq} have a strong solution $\bQ$ of the form
\eqref{QQ}.  Assume, in addition, that
\begin{enumerate}
\item[(i)] $\bQ$ is a Hilbert-Schmidt operator,

\item[(ii)] $ \bB\bQ (\bA-z)^{-1}$ is a trace class operator for
$\Img(z)\ne 0$,

\item[(iii)] $\bH$ and $\bA$ are resolvent comparable.
\end{enumerate}
Then   $A_i+B_{i j}Q_{ji}$,
$i=0,1, $ $j=1-i,$ are  admissible operators.
Moreover,
$A_i+B_{i j}Q_{ji}$ and $A_i$, $i=0,1, $ $j=1-i,$ are resolvent comparable.
 For  the spectral
shift function $\xi(\lambda, \bH, \bA)$ associated with the pair
of self-adjoint operators $(\bH,\bA)$ we have the decomposition
\begin{equation}
\label{ksi0}
\xi(\lambda;\bH, \bA)
 =\xi(\lambda; A_0+B_{01 }Q_{10},A_0)+
 \xi(\lambda; A_1+B_{10 }Q_{01},A_1)
\quad (\modulo\bbZ), \quad
\end{equation}
$$
\text{ for a.\,e. } \lambda \in \bbR.
$$
In particular, the operator matrix $\bH$ can be block
diagonalized by a unitary transformation
\eqref{decom}
$$
\bU^*\bH\bU=
\left(\begin{array}{cc}
H_0      &   0\\
0 &   H_1
\end{array}\right),
$$
where $\bU$ is the unitary operator from the polar decomposition
$$
\bI+\bQ=\bU|\bI+\bQ|,
$$
and
\begin{equation}
\label{ksi00}
\xi(\lambda;\bH, \bA)
 =\xi(\lambda; H_0,A_0)+
 \xi(\lambda; H_1,A_1)
\quad (\modulo\bbZ), \quad
\end{equation}
$$
\text{ for a.\,e. } \lambda \in \bbR.
$$
\end{theorem}
\begin{proof}
By Theorem \Ref{th42} (i) the normal operator $\bV=\bI+\bQ$
has a bounded inverse. Due to the
assumption (i) the spectrum of $\bV$ is purely point. Thus, by
Lemma\,\Ref{dixi}\,(vi)
\begin{equation}
\label{class}
\bV\in\cA(\cH),
\end{equation}
where $\cA(\cH)$ is the zero trace commutator class introduced by
Definition  \Ref{DefTC}.
By Theorem \Ref{thHi2} (ii) one concludes
\begin{equation}
\label{iso}
\bV^{-1}\bH\bV=\bA+\bB\bQ.
\end{equation}
Therefore, since by hypothesis (ii) the operator $\bB\bQ$ is a
relatively trace class perturbation of $\bA$,
one concludes that  the operators $\bV^{-1}\bH\bV$ and $\bA$ are
resolvent comparable. By condition (iii) $\bH$ and $\bA$ are
also resolvent comparable, and, therefore, by \eqref{class} the operator
$\bV^{-1}\bH\bV$ is admissible with the self-adjoint
representative $\bH$. Thus,  the stability property holds
\begin{equation}
\label{ksi1}
\xi(\lambda;\bV^{-1}\bH\bV, \bA)=\xi(\lambda;\bH, \bA)
\quad (\modulo\bbZ), \quad
\text{ for a.\,e. } \lambda \in \bbR,
\end{equation}
by the definition of the spectral shift function for resolvent comparable
admissible operators.

Next, let $\bV=\bU|\bV|$ be the polar decomposition of $\bV$.
By Theorem \Ref{thHi2} (iii) the operator $\bV\bV^*$ is diagonal
with respect to the decomposition \eqref{decom}. Using
representation \eqref{bV2} one infers that $\bV\bV^*-\bI$ is a
trace class operator, since $\bQ$ is the Hilbert--Schmidt
operator by the hypothesis.  Therefore,
\begin{equation}
\label{vib}
|\bV|-\bI\in \cB_1(\cH),
\end{equation}
where $|\bV|=(\bV\bV^*)^{1/2}$, and, hence,
$|\bV|\in \cA(\cH)$ by \eqref{vib}.  The operator ${\bf V^{-1}HV}$
is similar to the self-adjoint operator
$\bU^*\bH\bU$:
\begin{equation}
\label{posd}
{\bf V^{-1}HV}=|\bV|(\bU^*\bH\bU)|\bV|^{-1}.
\end{equation}
Using \eqref{vib} and \eqref{posd}, one concludes that
$\bV^{-1}\bH\bV$ and $\bU^*\bH\bU$ are resolvent comparable.
Therefore, taking into account that $|\bV|\in \cA(\cH)$ one infers that
 $\bU^*\bH\bU$ is a self-adjoint representative of the
admissible operator
 $\bV^{-1}\bH\bV$ and, hence,
\begin{equation}
\label{ksi2}
\xi(\lambda;\bV^{-1}\bH\bV, \bA)=\xi(\lambda;\bU^*\bH\bU, \bA)
\quad (\modulo\bbZ), \quad
\text{ for a.\,e. } \lambda \in \bbR,
\end{equation}
by Lemma \eqref{stability}.  By
Theorem \Ref{thHi2} (iii) the operator  $\bU^*\bH\bU$ is
diagonal with respect to decomposition \eqref{decom}
$$
\bU^*\bH\bU=\left(\begin{array}{cc}
H_0      &   0\\
0 &   H_1
\end{array}\right),
$$
where $H_i$, $ i=0,1,$ are self-adjoint operators in the Hilbert
spaces $\cH_i, $ $i=0,1,$ introduced by \eqref{Hip} and
\eqref{Hi}.  Since $\bU^*\bH\bU$ is a block-diagonal operator,
by additivity of the spectral shift function associated with a
pair of self-adjoint operators with respect to direct sum
decompositions (which follows from the definition of the
spectral shift function by the  trace formula \eqref{5.1}) one
obtains that
\begin{equation}
\label{ksi3}
\xi(\lambda; \bU^*\bV\bU, \bA)=\sum_{i=0}^1\xi(\lambda; H_i, A_i).
\end{equation}

By Theorem \Ref{thHi2} (ii) the operator  $\bV^{-1}\bH\bV$ is
diagonal with respect to the decomposition \eqref{decom}
$$
\bV^{-1}\bH\bV=\left(\begin{array}{cc}
A_0+B_{01}Q_{10}      &   0\\
0 &   A_1+B_{10}Q_{01}
\end{array}\right),
$$
where $A_i+B_{ij}Q_{ji}$, $ i=0,1$, $j=1-i,$ are operators
similar to self-adjoint operators $H_i$ given by \eqref{Hip}:
$$
H_i=(I_{\cH_i}+Q^*_{ji}Q_{ji})^{1/2}(A_i+B_{ij}Q_{ji})
(I_{\cH_i}+Q^*_{ji}Q_{ji})^{-1/2},\quad i=0,1,\quad
j=1-i.
$$
Here $Q_{ij}$, $i=0,1$, $j=1-i,$ are the entries in the matrix
representation for the operator $\bQ$
$$
\bQ=
\left(\begin{array}{cc}
  0    &   Q_{01}\\
Q_{10} &   0
\end{array}\right).
$$
By hypothesis (i)
$\bQ$ is a Hilbert-Schmidt operator, which proves that
\begin{align*}
(I_{\cH_i}+Q^*_{ji}&Q_{ji})^{1/2}-I_{\cH_i}
\in \cB_1(\cH_i)\\
& i=0,1, \quad j=1-i.
\end{align*}
Therefore, the operators $A_i+B_{ij}Q_{ji}$, $i=0,1$,
$j=1-i,$ are admissible with  the self-adjoint representatives $H_i$.
Since $\bV^{-1}\bH\bV$ and $\bA $  are resolvent
comparable, so $A_i+B_{ij}Q_{ji}$  and  $A_i$,  $i=0,1$,
$j=1-i,$ are.
Hence, we have the following representation by Lemma \Ref{stability}
\begin{equation}
\label{ksi4}
\xi(\lambda;H_i, A_i)=\xi(\lambda; A_i+B_{ij}Q_{ji} , A_i),
\quad (\modulo\bbZ), \quad
\text{ for a.\,e. } \lambda \in \bbR,
\end{equation}
$$
\quad i=0,1,\,\, j=1-i.
$$
Combining \eqref{ksi1},\eqref{ksi2}, \eqref{ksi3}, and
\eqref{ksi4} proves \eqref{ksi0}.
\end{proof}
\begin{remark}
If the operator $\bQ$ is a trace class operator,  the conditions
{\rm(ii)} and {\rm(iii)} hold automatically. Therefore, they are
redundant in this case.
\end{remark}

\section{Further properties of the spectral shift function}
\label{SecSuff}

Throughout this section we assume that the spectra of the main
diagonal entries $A_0$ and $A_1$ of the operator matrix
\eqref{twochannel} are separated. More specifically, we will
adopt one of the three following hypotheses.
\begin{hypothesis}
\label{HEnorm}
Assume Hypothesis \Ref{diag} and suppose that the separation condition
\begin{equation}
\label{dist}
\dist\{\spec(A_0),\spec(A_1)\}=d >0
\end{equation}
holds true.  Assume, in addition, that  $B_{10}$ has a finite
norm  with respect to the spectral measure of $A_0$ or/and
$A_1$ and, moreover,
\begin{equation}
\label{Vard}
\|B_{01}\|
  \min \{
\|B_{01}\|_{E_{A_1}},\|B_{01}\|_{E_{A_0}}\}
<\frac{d^2}{4}.
\end{equation}
\end{hypothesis}
\begin{hypothesis}
\label{HBpi}
Assume Hypothesis \Ref{diag} and suppose that the separation condition
\eqref{dist} holds true.  Assume, in addition, that both
operators $A_0$ and $A_1$ are bounded and
\begin{equation}
\label{Bpi}
  \|B_{01}\|<\frac{d}{\pi}.
\end{equation}
\end{hypothesis}
\begin{hypothesis}
\label{HAdL}
Assume Hypothesis \Ref{diag}. Assume, in addition, that
the operator $A_0$ is semibounded from above,
$$
A_0\leq a_0<+\infty,
$$
the operator $A_1$ is semibounded from below,
$$
A_1\geq a_1>-\infty,
$$
and
$$
a_0<a_1.
$$
\end{hypothesis}
\begin{theorem}
\label{Qsolv}
Assume Hypothesis \Ref{HEnorm}.
Then the block operator matrix $\bH$ has two {\rm(}orthogonal to each
other{\rm)} reducing graph subspaces $\cG_i=\cG(\cH_i, Q_{ji})$,
$i=0,1,$ $j=1-i,$ associated with angular
operators $Q_{ji}\in \cB(\cH_i,\cH_j)$ such that
$$
Q_{10}=-Q_{01}^*
$$
and
\begin{align}
\label{Star}
\|B_{ij}Q_{ji}\|&\leq\frac{d}{2}-\sqrt{\frac{d^2}{4}-
\|B\|_{01}{\rm min}\{\|B_{01}\|_{E_{A_0}},\|B_{01}\|_{E_{A_1}}\}}
%
<\frac{d}{2},\\
\nonumber
&\quad i=0,1,\quad j=1-i.
\end{align}
Moreover,
the graph subspaces $\cG_i$, $i=0,1,$ are the spectral subspaces
of $\bH$ and $\cG_0\oplus\cG_1=\cH.$
\end{theorem}
\begin{proof}
Assume, for definiteness, that the operator $B_{10}$ has a finite
norm with respect to the spectral measure of the diagonal entry
$A_1$  of $\bH$ and the inequality holds
\begin{equation}\label{defin}
\|B_{01}\| \, \|B_{10}\|_{E_{A_1}}<\frac{d^2}{4}.
\end{equation}
Recall that by definition
$\|B_{10}\|_{E_{A_1}}=\|B_{10}^*\|_{E_{A_1}}$ and hence
$\|B_{10}\|_{E_{A_1}}=\|B_{01}\|_{E_{A_1}}$.

By Theorem \Ref{QsolvN} (ii) the Riccati equation
\begin{equation}
\label{pervoe}
QA_0-A_1Q+QB_{01}Q=B_{10}
\end{equation}
has a unique strong solution $Q\in\cB(\cH_0, \cH_1).$ Therefore,
the dual Riccati equation
\begin{equation}
\label{vtoroe}
KA_1-A_0K+KB_{10}K=B_{01}
\end{equation}
has a unique strong solution $K\in\cB(\cH_1)$ by Theorem
\Ref{RicEC}, and, moreover, $K=-Q^*$. Introducing the notations
$Q_{10}=Q$ and $Q_{01}=K$, equations \eqref{pervoe} and
\eqref{vtoroe} can be rewritten in the form
\begin{equation}
\label{iundeksu}
Q_{ji}A_i-A_jQ_{ji}+Q_{ji}B_{ij}Q_{ji}=B_{ji}, \quad i=0,1,\,
j=1-i.
\end{equation}
Therefore, the Riccati equation \eqref{riceq} has a strong
solution of the form \eqref{QQ}. Applying Lemma \Ref{princip}
one proves that the subspaces $\cG_i=\cG(\cH_i, Q_{ji})$,
$i=0,1$ $j=1-i,$ are reducing subspaces for $\bH$, which proves
the first assertion of the theorem under hypothesis \eqref{defin}.

In the case where $B_{10}$ has a finite norm with
respect to the spectral measure of  the diagonal entry$A_0$ and the inequality
$$
\|B_{01}\| \, \|B_{10}\|_{E_{A_0}}<\frac{d^2}{4}
$$
holds, the proof can be performed in an analogous way.

 Applying Theorem  \ref{QsolvN} (ii) (eq. \eqref{QEst1})
 proves estimate \eqref{Star} which, in turn, proves that
$$
\dist\{\spec(A_0+B_{01}Q_{10}),\spec(A_1+B_{10}Q_{01})\}>0.
$$
The last assertion of the theorem is a corollary of Theorem \Ref{thHi2}.

The proof is complete.
\end{proof}
\begin{remark}
Under Hypothesis \Ref{HEnorm},
if
$$
\Sup\limits_{\mu\in\,\spec(A_j)}
\|\mu\,(A_i+B_{ij}Q_{ji}-\mu)^{-1}\|<\infty
$$
for some $i=0,1,$ $ j=1-i$, then the strong solutions of the
Riccati equations \eqref{pervoe}, \eqref{vtoroe} turn out to be
the operator solutions by Theorem \Ref{RicEC}.
\end{remark}
Under Hypothesis \Ref{HBpi} one has a similar result.

\begin{theorem}
\label{QsolvPi}
Assume Hypothesis \Ref{HBpi}.
Then the block operator matrix $\bH$ has two {\rm(}orthogonal to each
other{\rm)} reducing graph subspaces $\cG_i=\cG(\cH_i, Q_{ji})$,
$i=0,1,$ $j=1-i,$ associated with the strictly contractive angular
operators $Q_{ji}\in \cB(\cH_i,\cH_j)$, $\|Q_{ji}\|<1$, such that
$$
Q_{10}=-Q_{01}^*.
$$
Moreover, the graph subspaces $\cG_i$, $i=0,1,$ are the spectral subspaces
of $\bH$ and $\cG_0\oplus\cG_1=\cH.$
\end{theorem}
\begin{proof}
The proof is analogous to that of Theorem \Ref{Qsolv}. The only
difference is that now we refer to  part (i) of Theorem
\Ref{QsolvN}, since for bounded $A_i\in\cB(\cH_i)$, $i=0,1$.
the concepts of the  weak, strong and operator solutions of  the
Riccati equations \eqref{iundeksu} coincide.
\end{proof}

The following statement has been proven in \cite{AdLMSr}.
\begin{theorem}
\label{TAdL}
Assume Hypothesis \Ref{HAdL}. Then for any
$B_{01}\in\cB(\cH_1,\cH_0)$ and $B_{10}=B_{01}^*$ the open
interval $(a_0,a_1)$ appears to be a spectral gap for $\bH$.
At the same time the spectral subspaces of the operator $\bH$
corresponding to the intervals $(-\infty,a_0]$ and
$[a_1,+\infty)$ admit representation in the form of graph
subspaces associated with the pairs $(\cH_0,Q_{10})$ and
$(\cH_1,Q_{01})$ for some
$Q_{ij}\in\cB(\cH_j,\cH_i)$, $i=0,1,$ $j=1-i.$ That is,
\begin{equation}
\label{E0gs}
\ran\bigg(E_{\bH}\big((-\infty,a_0]\big)\bigg)=\cG(\cH_0, Q_{10})
\end{equation}
and
\begin{equation}
\label{E1gs}
\ran\bigg(E_{\bH}\big([a_1,+\infty)\big)\bigg)=\cG(\cH_1,Q_{01}),
\end{equation}
where $E_{\bH}(\Delta)$ denotes the spectral projection of $\bH$
associated with the Borel set $\Delta\subset \bbR$. The angular
operators $Q_{ij}$  are strict contractions, $\|Q_{ij}\|<1$,
possessing the property $Q_{10}=-Q_{10}^*$.

Moreover, the projections $E_{\bH}\big((-\infty,a_0]\big)$ and
$E_{\bH}\big([a_1,+\infty)\big)$ can be expressed in terms of the
operator $Q=Q_{01}=Q_{10}^*$ in the following way
$$
E_{\bH}\big ((-\infty,a_0]\big )
=\left(\begin{array}{cc}
  (I_0+QQ^*)^{-1}       &  -(I_0+QQ^*)^{-1}Q  \\
 -Q^*(I_0+QQ^*)^{-1}       &   Q^*(I_0+QQ^*)^{-1}Q
\end{array}\right)
$$
and
$$
E_{\bH}\big ([a_1,+\infty)\big )
=\left(\begin{array}{cc}
 Q(I_1+Q^*Q)^{-1}Q^*       &  Q(I_1+Q^*Q)^{-1}  \\
 (I_1+Q^*Q)^{-1}Q^*        &   (I_1+Q^*Q)^{-1}
\end{array}\right).
$$
\end{theorem}
\begin{corollary}
\label{CAdL-AMM}
Assume Hypothesis \Ref{HAdL}. Then for any
$B_{01}\in\cB(\cH_1,\cH_0)$ and $B_{10}=B_{01}^*$
the Riccati equation
\begin{equation}
\label{RicSym}
Q_{10}A_0-A_1 Q_{10}+Q_{10}B_{01}Q_{10}=B_{10}
\end{equation}
has a strong contractive solution $Q_{10}\in\cB(\cH_0,\cH_1)$,
$\|Q_{10}\|<1$, and $Q_{01}=-Q_{10}^*$ is the strong solution to
the dual Riccati equation
\begin{equation}
\label{RicSymK}
Q_{01}A_1-A_0 Q_{01}+Q_{01}B_{10}Q_{01}=B_{01}.
\end{equation}
For the spectra of the operators $A_0+B_{01}Q_{10}$ with
$\dom(A_0+B_{01}Q_{10})=\dom(A_0)$
and $A_1+B_{10}Q_{01}$ with $\dom(A_1+B_{10}Q_{01})=\dom(A_1)$
the following inclusions hold true:
\begin{equation}
\label{separ}
\spec(A_0+B_{01}Q_{10})\subset(-\infty,a_0]
\quad\text{and}\quad
\spec(A_1+B_{10}Q_{01})\subset[a_1,+\infty).
\end{equation}
\end{corollary}
\begin{proof}
Any spectral subspace for $\bH$ is its reducing subspace. Thus,
by Theorem \Ref{TAdL} the subspaces \eqref{E0gs} and \eqref{E1gs}
are reducing graph subspaces for $\cH$. Then Lemma \Ref{princip}
implies that the angular operators $Q_{01}$ and $Q_{10}$ from
the r.\,h.\,s. parts of formulas \eqref{E0gs} and \eqref{E1gs}
are strong solutions to equations  \eqref{RicSym} and
\eqref{RicSymK}, respectively. A proof of \eqref{separ} can be
found in \cite{AdLMSr}.
\end{proof}
\begin{remark}
Under Hypothesis \Ref{HAdL} the case where one of the self-adjoint
operators $A_0$ or $A_1$ is bounded has been treated first in
\cite{AdL}. Recently this case has been revisited in \cite{AdLT}
where sufficient conditions implying uniqueness of the strictly
contractive solutions to the operator Riccati equations have
been found.
\end{remark}
\begin{lemma}
\label{Qsolvc}
Assume at least one of the Hypotheses \Ref{HEnorm}, \Ref{HBpi},
and \Ref{HAdL}. Then the block operator matrix
$$
\bH_t=\bA+t\bB, \quad t\in [0,1]
$$
has two {\rm(}orthogonal to each other{\rm)}
reducing graph subspaces
$$
\cG(\cH_i,\sQ_{ji}(t)), \quad i=0,1,\,j=1-i,\,\, t\in[0,1],
$$
associated with angular operators
$\sQ_{ji}(t)\in \cB(\cH_i,\cH_j)$ which continuously depend on
$t\in [0,1]$ in the norm of the space $\cB(\cH_i,\cH_j)$.
In addition, under Hypothesis \Ref{HEnorm} the following holds true:
\begin{align*}
\|B_{ij}\sQ_{ji}(t)\|&\leq\frac{d}{2}-\sqrt{\frac{d^2}{4}-
\|B\|_{01}{\rm min}\{\|B_{01}\|_{E_{A_0}},\|B_{01}\|_{E_{A_1}}\}}
%
<\frac{d}{2},\\
&\quad i=0,1,\,\,j=1-i,\,\,t\in[0,1].
\end{align*}
Under Hypotheses \Ref{HBpi} or \Ref{HAdL} the operators $\sQ_{ji}(t)$
are strict contractions,
$$
\|\sQ_{ji}(t)\|<1, \quad
i=0,1, \quad j=1-i,\quad t\in[0,1].
$$
\end{lemma}
\begin{proof}
Under Hypothesis \Ref{HEnorm} or \Ref{HBpi} this assertion is an
immediate consequence of Theorems \Ref{Qsolv} or \Ref{QsolvPi}
respectively, and Remark \Ref{conti}.

Therefore, assume Hypothesis \Ref{HAdL}. Since the operator $\bB$
is bounded, and the interval $(a_0,a_1)$ does not contain points
of the spectrum of $\bH_t$ for all $t\in \bbR$, by a result by
Heinz \cite{H51} (see also \cite{Kato}, Theorem 5.12) the
spectral projection
$$
E(t)=E_{\bH_t}\big((-\infty,a_0]\big ), \quad t\in \bbR
$$
continuously depends on $t\in \bbR$ in the uniform operator
topology.  By Theorem \Ref{TAdL} the projection $E(t)$ admits
matrix representation with respect to the direct sum of the
Hilbert spaces $\cH_0\oplus\cH_1$
$$
E(t)=
\left(\begin{array}{cc}
  (I_0+Q_tQ^*_t)^{-1}       &  -(I_0+Q_tQ_t^*)^{-1}Q_t  \\
 -Q_t^*(I_0+Q_tQ^*_t)^{-1}       &   Q_t^*(I_0+Q_tQ_t^*)^{-1}Q_t
\end{array}\right),\quad t\in \bbR,
$$
where $Q_t=\sQ_{01}(t)$, $t\in \bbR$. In particular, the
continuity of the family $\{ E(t)\}_{t\in \bbR}$ implies the
continuity of the families of operators
$\{(I_0+Q_tQ_t^*)^{-1}\}_{t\in \bbR}$ and
$\{(I_0+Q_tQ_t^*)^{-1}Q_t\}_{t\in \bbR}$ in the uniform operator
topology of the spaces $\cB(\cH_0)$ and $\cB(\cH_1,\cH_0)$,
respectively.  Since the family $\{(I_0+Q_tQ_t^*)^{-1}\}_{t\in
\bbR}$ is continuous, the family $\{(I_0+Q_tQ_t^*)\}_{t\in
\bbR}$ is also continuous.  Multiplying the operator
$(I_0+Q_tQ_t^*)^{-1}Q_t$  by $I_0+Q_tQ_t^*$  from the left
proves the continuity of the angular operators $Q_t$ as a
function of $t$ in the uniform operator topology. Recalling now
that $\sQ_{10}(t)=-\sQ_{01}(t)^*=-Q_t^*$ proves the continuity
of the family $\sQ_{ij}(t)$, $i=0,1$, $j=1-i,$ as a function of
the parameter $t\in\bbR$ in the uniform operator topology. The
proof is complete.
\end{proof}

To a large extent, the angular operator $\bQ$, being a strong
solution to the Riccati equation \eqref{riceq}, inherits some
properties of the operator $\bB$.  For instance, if $\bB$
belongs to a symmetric ideal, so does $\bQ$, provided that the
certain spectra separation conditions are fulfilled for $A_0$
and $A_1$.  In fact, we have the following result (for
simplicity, formulated using the scale of Schatten--von Neumann
ideals).
\begin{theorem}
\label{pogonia}
Assume Hypothesis \Ref{diag} and let the Riccati equation
\eqref{riceq} have a strong solution $\bQ$ of the form
\eqref{QQ} with respect to the decomposition
$\cH=\cH_0\oplus\cH_1$.  Assume, in addition, that  either
condition \eqref{dist} is valid or the condition
\begin{equation}
\label{AH}
\dist\{\spec(A_i+B_{ij}Q_{ji}),\spec(A_j)\}> 0
\quad \text{ for some } i,j=0,1, \quad i\ne j,
\end{equation}
holds.
Then if $\bB\in \cB_p(\cH)$ for some $p\ge 1$,
then $\bQ\in \cB_p(\cH)$.
\end{theorem}
\begin{proof}
We recall that the strong solvability of the Riccati equation
\eqref{riceq} under constraint \eqref{QQ} is equivalent to the
strong solvability of the following pair of equations
\begin{equation}
\label{QAAQ}
Q_{ji}A_i-A_i Q_{ji}=B_{ji}-Q_{ji}B_{ij}Q_{ji}, \quad i=0,1, \,\,j=1-i.
\end{equation}
Therefore, the assumption $\bB\in \cB_p(\cH)$ for some $p\ge 1$
implies $B_{ij}\in \cB_p(\cH_j,\cH_i)$ $ i=0,1,$  $j=1-i$.
Hence, the r.\,h.\,s. of \eqref{QAAQ} is an element of the
space $\cB_p(\cH_i,\cH_j)$.  Under hypothesis \eqref{dist} one
concludes that $Q_{ji}\in\cB_p(\cH_i,\cH_j)$ by Theorem
\Ref{Kiexp} (in particular by estimate \eqref{ozenka}), and,
thus, $\bQ\in \cB_p(\cH)$, since \eqref{QQ} holds.

Further, assume that \eqref{AH} holds for some $i=0,1,$
$j=1-i$.  By Theorem \eqref{thHi2} the operator
$A_i+B_{ij}Q_{ji}$, $i=0,1,$
$j=1-i,$ is similar to a self-adjoint operator
$H_i$. That is, the representation holds
\begin{equation}
\label{simsim}
A_i+B_{ij}Q_{ji}=V H_i V^{-1}, \quad i=0,1, \,\,j=1-i,
\end{equation}
for some $V_i\in\cB(\cH_i)$ such that $V_i^{-1}\in \cB(\cH_i)$
(see \eqref{Hip}). Therefore,  \eqref{QAAQ}
can be rewritten in the form
$$
Q_{ji}V_i H_i V_i^{-1}-A_j Q_{ji}=B_{ji}
$$
and, hence, the operator $X_{ji}=Q_{ji}V_i$ is a strong solution
to the Sylvester equation
$$
X_{ji}H_i-A_j X_{ji}=B_{ji}V_i,\quad i=0,1, \,\,j=1-i.
$$
By \eqref{AH} and \eqref{simsim}
one infers
$$
\dist\{\spec(A_0+B_{01}Q_{10}),\spec(A_1)\}> 0.
$$
Meanwhile,  the assumption $\bB\in \cB_p(\cH)$ for some $p\ge 1$
implies $B_{ji}\in\cB_p(\cH_i,\cH_j)$ and, hence, $B_{ji}V_i\in
\cB_p(\cH_i,\cH_j)$, $i=0,1,$ $j=1-i$.  Applying Theorem
\Ref{Kiexp} once more, one deduces that
$X_{ji}\in\cB_p(\cH_i,\cH_j)$.  Hence,
$Q_{ji}=T_{ji}V_i^{-1}\in\cB_p(\cH_i,\cH_j)$, $i=0,1,$ $j=1-i$.
Finally, by \eqref{QQ} one concludes that $\bQ\in\cB_p(\cH)$.

The proof is complete.
\end{proof}

In what follows we need one abstract result of a topological
nature.
\begin{lemma}\label{abstract}
Let $L_t$, $t\in [0,1]$  be a one-parameter family of
self-adjoint operators such that $L_t$ and $L_0$ are resolvent
comparable for all $t\in [0,1]$ and the difference
$(L_t-z)^{-1}-(L_0-z)^{-1}$, $\Img (z)\ne 0$, is a continuous
function of $t\in [0,1]$ in the trace class topology.  Assume,
in addition, that
$$
[a,b]\cap\spec(L_t)=\emptyset \quad \text{ for all } t\in [0,1]
$$
for some $a,b \in \bbR,$ $a<b$.  Then for the unique family of
the spectral shift functions $\xi(\cdot;L_t, L_0)$ continuous in
$t\in [0,1]$ in the topology of the weighted space $L^1(\bbR;
(1+\lambda^2)^{-1})$ with the weight $(1+\lambda^2)^{-1}$ one
has
\begin{equation}
\label{vanish}
\xi(\lambda;L_t, L_0)=0 \text{ for a.\,e. } \lambda\in [a,b],\quad
t\in [0,1].
\end{equation}
\end{lemma}
\begin{proof}
The existence of the  one-parameter family of the spectral shift
functions $\xi(\cdot; L_t, L_0)$,  $t\in [0,1]$ that is
continuous in the topology of the weighted space $L^1(\bbR;
(1+\lambda^2)^{-1})$ is proven in \cite{Ya92}.  Next, since
$[a,b]$  belongs to the spectral gap of $L_t$ for any $t\in
[0,1]$, the spectral shift function $\xi(\lambda;L_t, L_0)$ is a
constant $n(t)\in \bbZ$ a.\,e. on the interval $[a,b]$. Integrating
the difference $n(t)-n(s)$ over $\lambda\in [a,b]$ with the
weight $(1+\lambda^2)^{-1}$ yields the estimate
$$
|n(t)-n(s)|\le \frac{\|\xi(\cdot\,;L_t,
L_0)-\xi(\cdot\,;L_s,
L_0)\|_{L^1(\bbR; (1+\lambda^2)^{-1})}}{\arctan(b)-\arctan(a)},
\quad t,s \in [0,1],
$$
which proves that $n(t) $ is a continuous integer-valued
function of $t\in [0,1]$. Since $n(0)=0$, it follows that
$n(t)=0$ for all $t\in [0,1]$.
\end{proof}

Now we are prepared to present the main result of the paper.
\begin{theorem}
\label{trsch}
Assume Hypothesis \Ref{diag} and at least  one  of
Hypotheses \Ref{HEnorm}, \Ref{HBpi}, and \Ref{HAdL} . Then the
Riccati equation \eqref{riceq} has a strong solution
$\bQ\in\cB(\cH)$ of the form
$$
\bQ=\left(\begin{array}{cc}
0      &   Q_{01}\\
Q_{10} &   0
\end{array}\right), \quad Q_{10}=-Q_{01}^*\in\cB(\cH_0,\cH_1),
$$
written with respect to the decomposition $\cH=\cH_0\oplus
\cH_1$ and hence the operator $\bH$ has reducing graph subspaces
$\cG_i=\cG(\cH_i, Q_{ji})$, $i=0,1,$ $j=1-i$.
If $\bH$ and $\bA$ are resolvent comparable and $\bB$ is a
Hilbert-Schmidt operator,
then   $A_i+B_{i j}Q_{ji}$,
$i=0,1, $ $j=1-i,$ are admissible operators.
Moreover,
$A_i+B_{i j}Q_{ji}$ and $A_i$, $i=0,1, $ $j=1-i,$ are resolvent comparable.
For  the spectral
shift function $\xi(\lambda, \bH, \bA)$ associated with the pair
of self-adjoint operators $(\bH,\bA)$ one has the decomposition
\begin{align}
\label{ksi000}
\xi(\lambda;\bH, \bA)
&=\xi(\lambda; A_0+B_{01 }Q_{10},A_0)+
\xi(\lambda; A_1+B_{10 }Q_{01},A_1)
\quad (\modulo\bbZ), \quad\\
\nonumber
&\text{ for a.\,e. } \lambda\in\bbR.
\end{align}
Moreover, the spectral shift functions
$\xi(\lambda;A_i+B_{ij}Q_{ji},A_i)$ associated with the pairs
$(A_i+B_{ij}Q_{ji},A_i)$, $i=0,1,$ $j=1-i,$ can be chosen in
such a way that
\begin{align}
\label{ksi001}
\xi(\lambda;A_i+B_{i j}Q_{ji}, A_i)&=0\quad
\text{ for a.\,e. }
\lambda\in \spec(A_j),\\
\nonumber
& i=0,1, \quad  j=1-i.
\end{align}
\end{theorem}
\begin{proof}
Under the assumptions of the theorem the existence of a strong
solution $\bQ\in\cB(\cH)$ of the Riccati equation \eqref{riceq}
is guaranteed by Lemma \Ref{princip} and Theorem \Ref{Qsolv},
Theorem \Ref{QsolvPi} or Corollary \Ref{CAdL-AMM}.
Since, by hypothesis, $\bB\in \cB_2(\cH)$, one infers $\bQ\in
\cB_2(\cH)$ by Theorem \Ref{pogonia}. Thus, the assumption (i)
of Theorem \Ref{th42} holds.  Therefore, $\bB\bQ$ is a trace
class operator, and hence the assumption (ii) of Theorem
\Ref{th42} holds. The assumption (iii) of Theorem \Ref{th42}
holds by hypothesis and, therefore, $A_i+B_{i j}Q_{ji}$, $i=0,1,
$ $j=1-i,$ are  admissible operators, $A_i+B_{i j}Q_{ji}$ and
$A_i$, $i=0,1, $ $j=1-i,$ are resolvent comparable and the
decomposition \eqref{ksi000} takes place by Theorem \Ref{th42}.

Introducing the family $\bH_t=\bA+t\bB$,  $t\in [0,1]$, by
Lemma \Ref{Qsolvc} one infers the existence of the
operators $\sQ_{ij}(t)\in \cB(\cH_i, \cH_j)$ that
continuously depend on $t\in[0,1]$ in the topology of the space
$\cB(\cH_i, \cH_j)$ and are such such
that $\bH_t$, $t\in[0,1]$ has reducing
graph subspaces
$$
\cG_i(t)=\cG(\cH_i, \sQ_{ji}(t)), \quad i=0,1,\,j=1-i,\,\, t\in
[0,1].
$$
Therefore, by Lemma \Ref{princip}
 the Riccati equation
\begin{equation}
\label{ricbold}
\bQ_t\bA-\bA\bQ_t+\bQ_t(t\bB)\bQ_t=t\bB,\quad t\in [0,1],
\end{equation}
has a strong solution $\bQ_t$ which reads
with respect to the decomposition~\eqref{decom} as
\begin{equation}
\label{qtt}
\bQ_t=\left(\begin{array}{cc}
0      &   \sQ_{01}(t)\\
\sQ_{10}(t) &   0
\end{array}\right),\quad t\in [ 0,1],
\end{equation}
and $\sQ_{ji}(t)=-[\sQ_{ij}(t)]^*$, $t\in
[0,1]$.  Hence, each entry $\sQ_{ji}(t)$, $t\in [0,1]$, in
\eqref{qtt} is a strong solution of the Riccati equation
\begin{equation}
\label{qt}
\sQ_{ji}(t) A_i-A_j \sQ_{ji}(t)=tB_{ji}-t\sQ_{ji}(t)B_{ij}\sQ_{ji}(t),
\quad t\in [0,1].
\end{equation}
Since $\sQ_{ji}(t)$ is continuous in the norm operator topology,
the r.\,h.\,s. of \eqref{qt} depends continuously on $t\in
[0,1]$ in the topology of the space $\cB_2(\cH_i,\cH_j)$.
Therefore, by Theorem \Ref{Kiexp} (estimate \eqref{ozenka}) the
path $\sQ_{ji}(t)$, $t\in [0,1]$, is continuous in the topology
of the space $\cB_2(\cH_j,\cH_i)$, and, thus, the family
$\{tB_{ij}\sQ_{ji}(t) \}_{t\in [0,1]}, \,$ $ i=0,1,$  $ j=1-i,$
is continuous in the topology of the space
$\cB_1(\cH_{i},\cH_j)$.

Clearly, the map
\begin{equation}
\label{map}
t\longrightarrow (A_i+tB_{ij}\sQ_{ji}(t)-z)^{-1}-
(A_i-z)^{-1})\in \cB_1(\cH_i), \quad t\in[0,1],
\end{equation}
$$
\quad i=0,1,\quad \Img(z)\ne 0.
$$
is continuous in the topology of the space $\cB_1(\cH_i)$,
$i=0,1.$ Taking into account that
 the family $\sQ^*_{ji}(t)\sQ_{ji}(t)$ is continuous in the
topology of $\cB_1(\cH_i)$, $i=0,1$, and
introducing the self-adjoint representatives of the admissible
operators $A_i+tB_{ij}\sQ_{ji}(t)$, $i=0,1, $ $ j=1-i$, $t\in
[0,1]$,
\begin{equation}
\label{h'(t)}
\sH_i(t)=[I_{\cH_i}+\sQ^*_{ji}(t)\sQ_{ji}(t)]^{1/2}(A_i+tB_{ij}\sQ_{ji}(t))
[I_{\cH_i}+\sQ^*_{ji}(t)\sQ_{ji}(t)]^{-1/2},
\end{equation}
$$
t\in[0,1],
$$
one concludes that the map
\begin{equation}\label{map1}
t\longrightarrow [\sH_i(t)-z)^{-1}- (A_i-z)^{-1}]\in \cB_1(\cH_i),
\quad t\in[0,1],
\end{equation}
is also continuous in the topology of  $\cB_1(\cH_i)$, $i=0,1.$

Let
$$
\Delta_i = \begin{cases}
\{\lambda:\,\dist\{\lambda,\spec(A_i)\}>d/2\},
& \text{if Hypothesis \Ref{HEnorm} holds,}\\
\{\lambda:\,\dist\{\lambda,\spec(A_i)\}>d/\pi\},
& \text{if Hypothesis \Ref{HBpi} holds,}\\
\bbR\,\backslash \, \overline{\text{
convex hull } (\spec(A_i))},
& \text{if Hypothesis \Ref{HAdL}
holds,}
\end{cases}
$$
$$
\quad i=0,1.
$$
Obviously
\begin{equation}
\label{spezii}
\spec(A_j)\subset \Delta_i, \quad i=0,1, \,\,j=1-i.
\end{equation}

Our claim is that $\Delta_i$, $i=0,1,$ belongs to the resolvent
set of $\sH_i(t)$,  $i=0,1,$ for all $t\in [0,1]$, that is,
\begin{equation}
\label{claim}
\Delta_i\cap \spec(\sH_i(t))=
\emptyset, \quad i=0,1, \,\,t\in
[0,1].
\end{equation}

Under Hypothesis \Ref{HAdL} the statement \eqref{claim} is a
consequence of Theorem \Ref{TAdL} (Eq. \eqref{separ}).

Assume, therefore, either  Hypotheses \Ref{HEnorm} or Hypotheses  \Ref{HBpi}.

Under Hypothesis \Ref{HEnorm}, applying Theorem \Ref{Qsolv}
one obtains the following uniform bounds
$$
\|tB_{ij}\sQ_{ji}(t)\|<\frac{d}{2},\quad t\in [0,1],\quad
i=0,1, \quad j=1-i.
$$
Thus, one concludes that
$$
\{\lambda:\,\dist\{\lambda,\spec(A_i)\}>d/2\}\bigcap\,
\spec(A_i+tB_{ij}\sQ_{ji}(t))=\emptyset
\quad \text{ for all } t\in [0,1],
$$
$$
i=0,1,\,\, j=1-i.
$$
Under Hypothesis \Ref{HBpi} the operator $\sQ_{ji}(t)$, $i=0,1$,
$j=1-i$, $t\in[0,1]$, is a strict contraction by Theorem
\Ref{QsolvPi}.  Therefore,
$$
\|tB_{ij}\sQ_{ji}(t)\|<\frac{d}{\pi},\quad t\in [0,1],\quad
i=0,1, \quad j=1-i,
$$
and
$$
\{\lambda:\,\dist\{\lambda,\spec(A_i)\}>d/\pi\}\bigcap\,
\spec(A_i+tB_{ij}\sQ_{ji}(t))=\emptyset
\quad \text{ for all } t\in [0,1],
$$
$$
i=0,1,\,\, j=1-i.
$$
By \eqref{h'(t)} the operators $\sH_i(t)$ and
$A_i+tB_{ij}\sQ_{ji}(t)$, $i=0,1,$ $j=1-i$, $t\in [0,1]$, are
similar to each other, which proves  \eqref{claim}  under
Hypotheses \Ref{HEnorm} or/and \Ref{HBpi}.

Applying Lemma \Ref{abstract} one proves that there is a family
of spectral shift functions $\xi(\,\cdot\,;\sH_i(t), A_i)\}_{t\in
[0,1]}$, $i=0,1,$ continuous in the topology of the weighted space
$L^2(\bbR;(1+\lambda^2)^{-1})$ such that
\begin{equation}
\label{nuliksi}
\xi(\lambda;\sH_i(t), A_i)=0 \text{ for a.\,e. } \lambda\in
[a_i,b_i],\quad t\in [0,1], \quad i=0,1,
\end{equation}
for any interval $[a_i,b_i]\subset \Delta_i$, $ i=0,1$. By
\eqref{map1} the operators  $(\sH_i(t)$ and $ A_i$, $i=0,1$,
$t\in [0,1]$, are resolvent comparable and, hence, by Lemma
\Ref{stability} one has the representation
$$
\xi(\lambda;A_i+tB_{ij}\sQ_{ji}(t), A_i)=
\xi(\lambda;\sH_i(t), A_i)\text{ for a.\,e. } \lambda\in \bbR,\quad
t\in [0,1], \quad i=0,1,
$$
since $\sH_i(t)$ are self-adjoint representatives of the admissible
operators  $A_i+tB_{ij}\sQ_{ji}(t)$, $i=0,1,$ $j=1-i$, $t\in [0,1]$.
It follows that the spectral shift functions $\xi(\lambda;
A_i+B_{ij}Q_{ji},A_i)$ associated with the pairs
$(A_i+B_{ij}Q_{ji}, A_i)$ $i=0,1, $ $j=1-i,$ can be chosen in
such a way that for any interval $[a_i,b_i]\subset\Delta_i,$ $i=0,1$,
\begin{equation}
\label{ksi01}
\xi(\lambda;A_i+B_{i j}Q_{ji}, A_i)=0\quad
\text{ for a.\,e. }\lambda\in [a_i,b_i]\subset \Delta_i
\end{equation}
$$
i=0,1, \quad  j=1-i,
$$
which, in particular, implies  assertion
\eqref{ksi001}, since \eqref{spezii} holds.
\end{proof}
\begin{remark}
\label{RemDet}
Assertion \eqref{ksi00} under Hypothesis \Ref{HAdL} in the
case where $\bB$ is a trace class operator has been proven by
Adamjan and Langer \cite{AdL00}.  Therefore, the main result of
the paper \cite{AdL00} in its part related to the existence of
the spectral shift function  and to the validity of
the representation \eqref{ksi00} is a particular case of our more
general considerations.
\end{remark}
\begin{corollary}
\label{final}
Assume the hypothesis of Theorem \Ref{trsch}. Then
\begin{enumerate}
\item[(i)]
the operator matrix $\bH$ can be block-diagonalized by a unitary
transformation
\eqref{decom}
$$
\bU^*\bH\bU=\left(\begin{array}{cc}
H_0      &   0\\
0 &   H_1
\end{array}\right),
$$
where $\bU$ is the unitary operator from the polar decomposition
$$\bI+\bQ=\bU|\bI+\bQ|\text{\rm ;}$$

\item[(ii)]  for the spectral shift function
$\xi(\lambda;\bH,\bA)$ the following splitting formula holds
\begin{align*}
\xi(\lambda;\bH, \bA) &=\xi(\lambda; H_0,A_0)+ \xi(\lambda;
H_1,A_1) \quad (\modulo\bbZ), \\
&\text{ for a.\,e. } \lambda \in \bbR\text{\rm ;}
\end{align*}

\item[(iii)] the spectral shift functions $\xi(\lambda;H_i,A_i)$,
$i=0,1,$ can be chosen in such a way that
\begin{equation}
\xi(\lambda; H_i, A_i)=0\quad
\text{ for a.\,e. }
\lambda\in \spec(A_{1-i}), \quad i=0,1.
\end{equation}
\end{enumerate}
\end{corollary}
\vspace*{2mm}
\noindent {\bf Acknowledgments.} A. K. Motovilov was supported
by the Deutsche For\-sch\-ungs\-gemeinschaft and by the Russian
Foundation for Basic Research.  He also gratefully acknowledges
the kind hospitality of the Institut f\"ur Angewandte
Mathematik, Universit\"at Bonn, during his stays in 2000 and
2001. K. A. Makarov is indebted to S. Fedorov, F.  Gesztesy, N.
Kalton, V. Kostrykin, and Yu. Latushkin for useful discussions.




\end{document}